\documentclass[A4paper,12pt]{book}
\usepackage{amssymb,amsmath}
\usepackage[russian]{babel}
\usepackage{graphicx}
\newcount\exnom\exnom=0

\long\def\exo/{\vspace{0.2cm} \noindent\advance\exnom by1{\bf
{\the\exnom}.}}
\newcommand{\z}{\exo/  }

\newcommand{\ds}{\displaystyle}

\def\im{\mathop{\rm Im}\nolimits}

\begin{document}

\thispagestyle{empty}

\noindent{\Large \bf \vspace{2mm}{Geometrical Structures of the Endomorphism Semiring\\ of a Finite Chain: Simplices, Strings and Triangles}}

\vspace{5mm}

\begin{quote}
{\bf \large Ivan Trendafilov}

\small
\emph{Department of Algebra and Geometry, Faculty of Applied Mathematics and\\ Informatics, Technical University of Sofia, Sofia, Bulgaria,\\ \emph{e-mail:}
ivan$\_$d$\_$trendafilov@abv.bg}

We establish new results concerning endomorphisms of a finite chain if the cardinality of the image of such endomorphism is no more than some fixed number $k$. The semiring of all such endomorphisms can be seen as a $k$--simplex whose vertices are the constant endomorphisms. We explore the properties of these $k$--simplices and find some results for arbitrary $k$. For $k=1$ and 2 we give a full description of simplices called strings and triangles, respectively.
\end{quote}

\vspace{5mm}

 \noindent{\bf \large 1. \hspace{0.5mm} Introduction and Preliminaries}

\vspace{3mm}

The endomorphism semiring of a finite semilattice is well studied in [1] -- [7].
In the present paper we give a new treatment of the subsemirings of  endomorphism semiring $\widehat{\mathcal{E}}_{\mathcal{C}_n}$ of a finite chain. We investigate  endomorphisms $\alpha \in \widehat{\mathcal{E}}_{\mathcal{C}_n}$ such that $|\im(\alpha)| \leq k$, where $k$ is a positive integer and $k \leq n$. The set of these endomorphisms is a $k$--simplex with vertices  constant endomorphisms $\overline{a} = \wr\, \underbrace{a, \ldots, a}_n\,\wr$ and proper sides all the $m$--simplices, where $m < k$, $m \neq 0$.

The paper is organized as follows. After the introduction and preliminaries, in  the second section we study $k$--simplices for any natural $k$. Although we do not speak about any distance here, we define discrete neighborhoods with respect to any vertex of the simplex. In the main result of the section Theorem 5, we prove that the biggest discrete neighborhoods of the least and biggest vertex of the simplex are the semirings of a special type. In the third section we start the inductive study of the simplices and consider  1-simplices called strings. Section 4 is devoted to the basic properties of   triangle $\triangle^{(n)}\{a,b,c\}$. Here we prove that any element of the interior of a triangle is a sum of the elements of two strings of this triangle.
In the next section we consider the set of endomorphisms of $\triangle^{(n)}\{a,b,c\}$ such that some of elements $a$, $b$ and $c$ occur just $m$ times, $0 \leq m < n$, of the image of the endomorphism. These sets are called layers of the triangle with respect to some vertex. When two boundary elements of the layer are idempotents, we call this layer a basic layer and prove that all basic layers with respect to $a$ and $c$ are semirings isomorphic to well-defined strings.

In the last section are placed the main results of the paper. Here we construct a new subsemiring of $\triangle^{(n)}\{a,b,c\}$, the so-called idempotent triangle, containing all the right identities of  triangle $\triangle^{(n)}\{a,b,c\}$. In Theorem 33 we prove that the idempotent triangle is a disjoint union of two subsmirings.  The first one is the semiring of the right identities of $\triangle^{(n)}\{a,b,c\}$ and the second is the maximal ideal of the idempotent triangle. In this section we describe the subsemirings of $a$--nilpotent, $b$--nilpotent and $c$--nilpotent endomorphisms as  geometric figures: trapezoids and parallelogram. In this section is also investigated two subsemirings which are geometric parallelograms. In Theorem 40 we show that any triangle  $\triangle^{(n)}\{a,b,c\}$ is a disjoint union of eight subsemirings.

\vspace{3mm}

 Since the terminology  for
semirings  is not completely standardized, we say what our conventions are.
 An algebra $R = (R,+,.)$ with two binary operations $+$ and $\cdot$ on $R$, is called a {\emph{semiring}} if:

$\bullet\; (R,+)$ is a commutative semigroup,

$\bullet\; (R,\cdot)$ is a semigroup,

$\bullet\;$ both distributive laws hold $ x\cdot(y + z) = x\cdot y + x\cdot z$ and $(x + y)\cdot z = x\cdot z + y\cdot z$ for any $x, y, z \in R$.

 Let $R = (R,+,.)$ be a semiring.
If a neutral element $0$ of the semigroup $(R,+)$ exists and $0x = 0$, or $x0 = 0$, it is called a {\emph{left}} or a {\emph{right zero}}, respectively, for all $x \in R$.
If $0\cdot x = x\cdot 0 = 0$ for all $x \in R$, then it is called {\emph{zero}}. An element $e$ of a semigroup $(R,\cdot)$ is called a {\emph{left (right) identity}} provided
that $ex = x$, or $xe = x$, respectively, for all  $x \in R$.
If a neutral element $1$ of the semigroup $(R,\cdot)$ exists, it is called {\emph{identity}}.

A nonempty subset $I$ of  $R$ is called an \emph{ideal} if  $I + I \subseteq I$, $R\,I \subseteq I$ and $I\, R \subseteq I$.

The facts concerning semirings can be found in [8]. For semilattices we refer to [9].

\vspace{2mm}

For a join-semilattice $\left(\mathcal{M},\vee\right)$  set $\mathcal{E}_\mathcal{M}$ of the endomorphisms of $\mathcal{M}$ is a semiring
 with respect to the addition and multiplication defined by:

 $\bullet \; h = f + g \; \mbox{when} \; h(x) = f(x)\vee g(x) \; \mbox{for all} \; x \in \mathcal{M}$,

 $\bullet \; h = f\cdot g \; \mbox{when} \; h(x) = f\left(g(x)\right) \; \mbox{for all} \; x \in \mathcal{M}$.

 This semiring is called the \emph{ endomorphism semiring} of $\mathcal{M}$.
In this article all semilattices are finite chains. Following [5] and [6], we fix a finite chain $\mathcal{C}_n = \; \left(\{0, 1, \ldots, n - 1\}\,,\,\vee\right)\;$
and denote the endomorphism semiring of this chain with $\widehat{\mathcal{E}}_{\mathcal{C}_n}$. We do not assume that $\alpha(0) = 0$ for arbitrary
$\alpha \in \widehat{\mathcal{E}}_{\mathcal{C}_n}$. So, there is not a zero in  endomorphism semiring $\widehat{\mathcal{E}}_{\mathcal{C}_n}$. Subsemirings
${\mathcal{E}}_{\mathcal{C}_n}^a$, where $a \in \mathcal{C}_n$, of $\widehat{\mathcal{E}}_{\mathcal{C}_n}$ consisting of all endomorphisms $\alpha$ with
property $\alpha(a) = a$ are considered in [4].

\vspace{2mm}

If $\alpha \in \widehat{\mathcal{E}}_{\mathcal{C}_n}$ such that $f(k) = i_k$ for any  $k \in \mathcal{C}_n$ we denote $\alpha$ as an ordered $n$--tuple
$\wr\,i_0,i_1,i_2, \ldots, i_{n-1}\,\wr$. Note that mappings will be composed accordingly, although we shall usually give preference to writing mappings on
the right, so that $\alpha \cdot \beta$ means ``first $\alpha$, then $\beta$''. The identity $\mathbf{i} = \wr\,0,1, \ldots, n-1\,\wr$ and all constant
endomorphisms $\kappa_i = \wr\, i, \ldots, i\,\wr$ are obviously (multiplicatively) idempotents.

Let $a \in \mathcal{C}_n$. For every  endomorphism $\overline{a} = \wr\,a \,a\, \ldots\, a\,\wr$ the elements of
$$\mathcal{N}_n^{\,[a]} = \{\alpha \; | \; \alpha \in \widehat{\mathcal{E}}_{\mathcal{C}_n}, \; \alpha^{n_a} = \overline{a} \; \mbox{for some natural number} \; n_a\}$$
are called $a$--\emph{nilpotent endomorphisms}. An important result for $a$--nilpotent endomorphisms is

\vspace{3mm}

\textbf{Theorem} (Theorem 3.3 from [3]), \textsl{For any natural $n$, $n \geq 2$, and $a \in \mathcal{C}_n$ the set of $a$-- nilpotent endomorphisms  $\mathcal{N}_n^{\,[a]}$ is a subsemiring of  $\widehat{\mathcal{E}}_{\mathcal{C}_n}$.  The order of this semiring is $\left|\mathcal{N}_n^{\,[k]}\right| = C_k . C_{n-k-1}$, where $C_k$ is the $k$ -- th Catalan number.}

\vspace{3mm}

Another useful result is

\vspace{3mm}

\textbf{Theorem} (Theorem 9 from [2]), \textsl{The subset of $\widehat{\mathcal{E}}_{\mathcal{C}_n}$, $n \geq 3$, of all idempotent endomorphisms with $s$ fixed points
$k_1, \ldots, k_s$, ${1 \leq s \leq n-1}$, is a semiring of order  $\ds \prod_{m=1}^{s-1} (k_{m+1} - k_{m})$.}

For definitions and results concerning simplices we refer to [10] and [11].

\vspace{5mm}

 \noindent{\bf \large 2. \hspace{0.5mm} Simplices}

\vspace{3mm}

Let us fix  elements $a_0, a_1, \ldots, a_{k-1} \in \mathcal{C}_n$, where $k \leq n$, $a_0 < a_1 < \ldots < a_{k-1}$, and let $A = \{a_0, a_1, \ldots, a_{k-1}\}$.
We consider  endomorphisms $\alpha \in \widehat{\mathcal{E}}_{\mathcal{C}_n}$ such that $\im(\alpha) \subseteq A$. Let us call   ${k}$ -- \emph{simplex} (or only a {\emph{simplex}) the set of all such endomorphisms $\alpha$. We denote  $k$ -- simplex by
 $\sigma^{(n)}_k(A) = \sigma^{(n)}\{a_0, a_1, \ldots, a_{k-1}\}$.

It is easy to see that $\im(\alpha) \subseteq A$ and $\im(\beta) \subseteq A$ imply $\im(\alpha+ \beta) \subseteq A$ and $\im(\alpha\cdot \beta) \subseteq A$. Hence, we find

\vspace{3mm}

\textbf{Proposition} \z  \textsl{For any $A = \{a_0, a_1, \ldots, a_{k-1}\} \subseteq \mathcal{C}_n$  $k$ -- simplex $\sigma^{(n)}_k(A)$ is a subsemiring of $\widehat{\mathcal{E}}_{\mathcal{C}_n}$.}

\vspace{3mm}

The number $k$ is called a {\emph{dimension}} of  $k$ -- simplex $\sigma^{(n)}_k(A)$. Endomorphisms $\overline{a_j}$, where $j = 0, \ldots, k - 1$, such that $\overline{a_j}(i) = a_j$ for any $i \in \mathcal{C}_n$ are called {\emph{vertices}} of  $k$ -- simplex $\sigma^{(n)}_k(A)$.

Any  simplex $\sigma^{(n)}\{b_0, b_1, \ldots, b_{\ell - 1}\}$, where $b_0, \ldots, b_{\ell-1} \in A$, is called a {\emph{face}} of  $k$ -- simplex $\sigma^{(n)}_k(A)$. If $\ell < k$,  face $\sigma^{(n)}\{b_0, b_1, \ldots, b_{\ell - 1}\}$ is called a {\emph{proper face}}.

\vspace{1mm}

 The proper faces of $k$ -- simplex $\sigma^{(n)}_k(A) = \sigma^{(n)}\{a_0, a_1, \ldots, a_{k-1}\}$ are:

$\bullet$ $0$ -- simplices, which are vertices $\overline{a_0}, \ldots, \overline{a_k}$.

$\bullet$ $1$ -- simplices, which are called {\emph{strings}}. They are denoted by $\mathcal{STR}^{(n)}\{a,b\}$, where $a, b \in A$.

$\bullet$ $2$ -- simplices, which are called {\emph{triangles}}. They are denoted by ${\triangle}^{(n)}\{a,b,c\}$, where $a, b, c \in A$.

$\bullet$ $3$ -- simplices, which are called {\emph{tetrahedra}}. They are denoted by $\mathcal{TETR}^{(n)}\{a,b,c,d\}$, where $a, b, c, d \in A$.

$\bullet$  The last proper  faces are $k-1$ -- simplices $\sigma^{(n)}_{k-1}(B)$, where $B = \{b_0, \ldots, b_{k-2}\} \subset A$.

\vspace{1mm}

 The {\emph{boundary}} of  $k$ -- simplex  $\sigma^{(n)}_k(A)$ is  a union of all its proper faces and is denoted by $\mathcal{BD}\left(\sigma^{(n)}_k(A)\right)$. The set $\mathcal{INT}\left(\sigma^{(n)}_k(A)\right) = \sigma^{(n)}_k(A) \backslash \mathcal{BD}\left(\sigma^{(n)}_k(A)\right)$ is called an {\emph{interior}} of  $k$ -- simplex  $\sigma^{(n)}_k(A)$. The boundary and the interior of $k$ -- simplex are, in general, not semirings.

\vspace{1mm}

For any natural $n$, endomorphism  semiring $\widehat{\mathcal{E}}_{\mathcal{C}_n}$ is  $n$ -- simplex with vertices $\overline{0}, \ldots, \overline{n-1}$. The interior of this simplex consists of endomorphisms $\alpha$, such that $\im(\alpha) = \mathcal{C}_n$. Since the latter is valid only for identity $\mathbf{i} = \wr\,0,1, \ldots n -1\, \wr$, it follows that ${\mathcal{INT}\left( \widehat{\mathcal{E}}_{\mathcal{C}_n}\right) = \mathbf{i}}$.
\vspace{2mm}

There is a partial ordering of the faces of dimension $k - 1$ of $k$ -- simplex by the following way:  least face does not contain the vertex $\overline{a_k}$ and  biggest face does not contain the vertex $\overline{a_0}$.

\vspace{2mm}

 The biggest face of the $n$ -- simplex $\widehat{\mathcal{E}}_{\mathcal{C}_n}$ is the $(n-1)$ -- simplex $\sigma^{(n)}_{n-1} \{1,  \ldots, n-1\}$.
Now $\widehat{\mathcal{E}}_{\mathcal{C}_n}\backslash \sigma^{(n)}_{n-1} \{1,  \ldots, n-1\} = {\mathcal{E}}_{\mathcal{C}_n}^{(0)}$ which is a subsemiring of $\widehat{\mathcal{E}}_{\mathcal{C}_n}$. Similarly, the least face of $\widehat{\mathcal{E}}_{\mathcal{C}_n}$ is $\sigma^{(n)}_{n-1} \{0,  \ldots, n-2\}$.
Then $\widehat{\mathcal{E}}_{\mathcal{C}_n}\backslash \sigma^{(n)}_{n-1} \{0,  \ldots, n-2\} = {\mathcal{E}}_{\mathcal{C}_n}^{(n-1)}$ which is also a subsemiring of $\widehat{\mathcal{E}}_{\mathcal{C}_n}$. The other faces of $\widehat{\mathcal{E}}_{\mathcal{C}_n}$, where $n \geq 3$, do not have this property. Indeed,
one middle side is $\sigma^{(n)}_{n-1} \{0,  \ldots, k-1, k+1, \ldots, n-1\}$. But set $R = \widehat{\mathcal{E}}_{\mathcal{C}_n}\backslash \sigma^{(n)}_{n-1} \{0,  \ldots, k-1, k+1, \ldots, n-1\}$ is not a semiring because for any $n \geq 3$ and any $k \in \{1, \ldots, n-2\}$ if $\alpha = \wr\, 0, \ldots, 0, k\,\wr \in R$, then $\alpha^2 = \overline{0} \notin R$.

Let us fix  vertex $\overline{a_m}$, where $m = 0, \ldots, k-1$, of   simplex $\sigma^{(n)}\{a_0, a_1, \ldots, a_{k-1}\}$.
The set of all endomorphisms $\alpha \in \sigma^{(n)}_k(A)$ such that $\alpha(i) = a_m$ just for $s$ elements $i \in \mathcal{C}_n$ is called {\emph{$s$-th layer of  $k$ -- simplex with respect to $\overline{a_m}$}}, where $s = 0, \ldots, n-1$. We denote the $s$-th layer of the  simplex with respect to $\overline{a_m}$ by $\mathcal{L}^{s}_{a_m}\left(\sigma^{(n)}\{a_0, a_1, \ldots, a_{k-1}\}\right)$.
 So, the $0$ -- layer with respect to any vertex of the $k$ -- simplex is a face of the $k$ -- simplex, hence, it is a semiring.  In the general case, the $s$-th layer $\mathcal{L}^{s}_{a_m}\left(\sigma^{(n)}\{a_0, a_1, \ldots, a_{k-1}\}\right)$, where $s \in \mathcal{C}_n$, $s = 1, \ldots, n-2$, is not a subsemiring of $k$ -- simplex.

\vspace{3mm}

From topological point of view,  set $\{\overline{a_m}\} \cup \mathcal{L}^{n-1}_{a_m}\left(\sigma^{(n)}\{a_0, a_1, \ldots, a_{k-1}\}\right)$ is a discrete neighborhood consisting of the ``nearest points to  point'' $\overline{a_m}$. We denote this set by $\mathcal{DN}^{\,1}_m$. Similarly, we define $\mathcal{DN}^{\,2}_m = \mathcal{DN}^{\,1}_m\cup \mathcal{L}^{n-2}_{a_m}\left(\sigma^{(n)}\{a_0, a_1, \ldots, a_{k-1}\}\right)$ and, in general,
$\displaystyle \mathcal{DN}^{\,t}_m = \{\overline{a_m}\}\cup \bigcup_{\ell = n-t}^{n-1} \mathcal{L}^{\ell}_{a_m}\left(\sigma^{(n)}\{a_0, a_1, \ldots, a_{k-1}\}\right)$, where $m = 0, \ldots, k-1$, $t = 1, \ldots, n$.

\vspace{3mm}

\textbf{Proposition} \z  \textsl{Let $\overline{a_m}$, where $m = 0, \ldots, k-1$, be a vertex of the  simplex $\sigma^{(n)}\{a_0, a_1, \ldots, a_{k-1}\}$ and $\mathcal{L}^{n-1}_{a_m}\left(\sigma^{(n)}\{a_0, a_1, \ldots, a_{k-1}\}\right)$ be the $(n-1)$-th layer of the $k$ -- simplex with respect to $\overline{a_m}$. Then the set $\mathcal{DN}^{\,1}_m = \{\overline{a_m}\} \cup \mathcal{L}^{n-1}_{a_m}\left(\sigma^{(n)}\{a_0, a_1, \ldots, a_{k-1}\}\right)$, where $m = 0, \ldots, k-1$, is a  subsemiring of $\sigma^{(n)}\{a_0, a_1, \ldots, a_{k-1}\}$.}

\vspace{1mm}

\emph{Proof.} We consider three cases.
\vspace{1mm}

\emph{Case 1.} Let $m = 0$. Then elements of $\mathcal{DN}^{\,1}_0$ are endomorphisms:
$$\overline{a_0}, \; (a_0)_{n-1}a_1 = \wr\, \underbrace{a_0, \ldots, a_0}_{n-1},a_1\, \wr,\; \ldots \,, \; (a_0)_{n-1}a_{k-1} = \wr\, \underbrace{a_0, \ldots, a_0}_{n-1},a_{k-1} \wr. $$
Since $\overline{a_0} < (a_0)_{n-1}a_1 < \cdots < (a_0)_{n-1}a_{k-1}$, it follows that  set $\mathcal{DN}^{\,1}_0$ is closed under the addition.

We find $(a_0)_{n-1}a_i \cdot \overline{a_0} = \overline{a_0} \cdot (a_0)_{n-1}a_i  = \overline{a_0}$ for all $i = 1, \ldots, k-1$. Also we have
$(a_0)_{n-1}a_i \cdot (a_0)_{n-1}a_j = (a_0)_{n-1}a_j \cdot (a_0)_{n-1}a_i = \overline{a_0}$ for all $i,j \in \{1, \ldots, k -1\}$ with the only exception when $a_{k-1} = n-1$. Now $\left((a_0)_{n-1}(n-1)\right)^2 = (a_0)_{n-1}(n-1)$, $(a_0)_{n-1}(n-1)\cdot (a_0)_{n-1}a_i = (a_0)_{n-1}a_i$ and $(a_0)_{n-1}a_i \cdot (a_0)_{n-1}(n-1) = \overline{a_0}$. Hence, $\mathcal{DN}^{\,1}_0$ is a semiring.

\vspace{1mm}

\emph{Case 2.} Let $m = k-1$. Then elements of $\mathcal{DN}^{\,1}_{k-1}$ are endomorphisms:
$$ a_0(a_{k-1})_{n-1} = \wr\,a_0, \underbrace{a_{k-1}, \ldots, a_{k-1}}_{n-1}\,\wr,\; \ldots \,, \; a_{k-2}(a_{k-1})_{n-1} = \wr\,a_{k-2}, \underbrace{a_{k-1}, \ldots, a_{k-1}}_{n-1}\,\wr,\; \overline{a_{k-1}}. $$
Since $a_0(a_{k-1})_{n-1} < \cdots < a_{k-2}(a_{k-1})_{n-1} < \overline{a_{k-1}}$, it follows that the set $\mathcal{DN}^{\,1}_{k-1}$ is closed under the addition.

We find $a_i(a_{k-1})_{n-1} \cdot \overline{a_{k-1}} = \overline{a_{k-1}}\cdot a_i(a_{k-1})_{n-1} = \overline{a_{k-1}}$ for all $i = 1, \ldots, k-1$.
 Also we have $a_i(a_{k-1})_{n-1}  \cdot a_j(a_{k-1})_{n-1}  = a_j(a_{k-1})_{n-1}  \cdot a_i(a_{k-1})_{n-1}  = \overline{a_{k-1}}$ for all ${i,j \in \{0, \ldots, k -2\}}$ with also the only exception when $a_0 = 0$. We have $\left(0(a_{k-1})_{n-1}\right)^2 = 0(a_{k-1})_{n-1}$, $0(a_{k-1})_{n-1}\cdot a_i(a_{k-1})_{n-1} = a_i(a_{k-1})_{n-1}$ and $a_i(a_{k-1})_{n-1} \cdot 0(a_{k-1})_{n-1} = \overline{a_{k-1}}$
Hence, $\mathcal{DN}^{\,1}_{k-1}$ is a semiring.

\vspace{1mm}

\emph{Case 3.}  Let $0 < m < k-1$. Then elements of $\mathcal{DN}^{\,1}_m$ are endomorphisms:
$$ a_0(a_m)_{n-1} = \wr\,a_0, \underbrace{a_m, \ldots, a_m}_{n-1}\,\wr,\; \ldots \,, \; a_{m-1}(a_m)_{n-1} = \wr\,a_{m-1}, \underbrace{a_m, \ldots, a_m}_{n-1}\,\wr,\; \overline{a_m}, $$
$$ (a_m)_{n-1}a_{m+1} = \wr\,\underbrace{a_m, \ldots, a_m}_{n-1},a_{m+1}\,\wr,\; \ldots \,, \; (a_m)_{n-1}a_{k-1} = \wr\,\underbrace{a_m, \ldots, a_m}_{n-1}, a_{k-1}\,\wr.$$
Since $a_0(a_m)_{n-1} < \cdots < a_{m-1}(a_m)_{n-1} < \overline{a_m} < (a_m)_{n-1}a_{m+1} < \cdots < (a_m)_{n-1}a_{k-1}$, it follows that  set $\mathcal{DN}^{\,1}_m$ is closed under the addition.

\vspace{1mm}

Now there are four possibilities:

\vspace{1mm}

\emph{3.1.} Let $0 < a_0$ and $a_{k-1} < n-1$. Then
$$a_i(a_m)_{n-1}\cdot a_j(a_m)_{n-1} = a_j(a_m)_{n-1}\cdot a_i(a_m)_{n-1} = \overline{a_m} \;\; \mbox{for any} \; i, j = 0, \ldots m-1,$$
$$(a_m)_{n-1}a_i \cdot (a_m)_{n-1}a_j = (a_m)_{n-1}a_j \cdot (a_m)_{n-1}a_i = \overline{a_m} \;\; \mbox{for any} \; i, j = m+1, \ldots k-1,$$
$$a_i(a_m)_{n-1}\cdot (a_m)_{n-1}a_j = (a_m)_{n-1}a_j\cdot a_i(a_m)_{n-1} = \overline{a_m}$$ $$ \mbox{for any} \; i = 0, \ldots, m-1 \; \mbox{and}\;  j = m+1, \ldots k-1.
$$

Since $a_i(a_m)_{n-1}\cdot \overline{a_m} = \overline{a_m}\cdot a_i(a_m)_{n-1} = \overline{a_m}$ for $i = 1, \ldots, m - 1$ and, in a  similar way, ${(a_m)_{n-1}a_j \cdot \overline{a_m} = \overline{a_m}\cdot (a_m)_{n-1}a_j = \overline{a_m}}$ for $j = m+1, \ldots, k -1$ and also $\left(\overline{a_m}\right)^2 = \overline{a_m}$, it follows that $\mathcal{DN}^{\,1}_m$ is a commutative semiring.
\vspace{1mm}

\emph{3.2.} Let $a_0 = 0$ and $a_{k-1} < n - 1$. Then $\left(0(a_m)_{n-1}\right)^2 = 0(a_m)_{n-1}$,
 $$ 0(a_m)_{n-1} \cdot a_i(a_m)_{n-1} = a_i(a_m)_{n-1}, \;a_i(a_m)_{n-1}\cdot 0(a_m)_{n-1} = \overline{a_m} \;\;\mbox{for any}\; i = 1, \ldots m-1\; \mbox{and}$$  $$0(a_m)_{n-1} \cdot (a_m)_{n-1}a_j = (a_m)_{n-1}a_j \cdot 0(a_m)_{n-1} =  \overline{a_m}\;\; \mbox{for any}\; j = m+ 1, \ldots, k -1.$$

 We also observe that $\overline{a_m}\cdot 0(a_m)_{n-1} = 0(a_m)_{n-1} \cdot \overline{a_m} = \overline{a_m}$. All the other equalities between the products of the elements of $\mathcal{DN}^{\,1}_m$ are the same as in \emph{3.1}.
\vspace{1mm}

 \emph{3.3.} Let $a_0 > 0$ and $a_{k-1} = n-1$. Then $\left((a_m)_{n-1}(n-1)\right)^2 = (a_m)_{n-1}(n-1)$,
 $$(a_m)_{n-1}(n-1) \cdot a_i(a_m)_{n-1} = a_i(a_m)_{n-1}\cdot (a_m)_{n-1}(n-1) = \overline{a_m}\;\;\mbox{for any}\; i = 1, \ldots m-1\; \mbox{and}$$
$$(a_m)_{n-1}(n-1) \cdot (a_m)_{n-1}a_j = (a_m)_{n-1}a_j, \; (a_m)_{n-1}a_j \cdot (a_m)_{n-1}(n-1) = \overline{a_m}$$
for any $j = m+ 1, \ldots, k -1$.

 We also observe that $\overline{a_m}\cdot (a_m)_{n-1}(n-1) = (a_m)_{n-1}(n-1) \cdot \overline{a_m} = \overline{a_m}$. All the other equalities between the products of the elements of $\mathcal{DN}^{\,1}_m$ are the same as in \emph{3.1}.
\vspace{1mm}

 \emph{3.4.} Let $a_0 = 0$ and $a_{k-1} = n-1$. Now all  equalities between the products of the elements of $\mathcal{DN}^{\,1}_m$ are the same as in \emph{3.1.}, \emph{3.2.} and \emph{3.3}. So,  $\mathcal{DN}^{\,1}_m$ is a semiring. \hfill $\Box$

 \vspace{3mm}

Any simplex $\sigma^{(n)}\{b_0, b_1, \ldots, b_{\ell - 1}\}$ which is a face of  simplex $\sigma^{(n)}\{a_0, a_1, \ldots, a_{k-1}\}$ is called
{\emph{internal of the simplex}} $\sigma^{(n)}\{a_0, a_1, \ldots, a_{k-1}\}$ if $a_0 \notin \sigma^{(n)}\{b_0, b_1, \ldots, b_{\ell - 1}\}$ and $a_{k-1} \notin \sigma^{(n)}\{b_0, b_1, \ldots, b_{\ell - 1}\}$. Similarly  simplex $\sigma^{(n)}\{a_0, a_1, \ldots, a_{k-1}\}$, which is a face of  $n$ -- simplex
$\widehat{\mathcal{E}}_{\mathcal{C}_n}$, is called {\emph{internal simplex}} if $0 \notin \sigma^{(n)}\{a_0, a_1, \ldots, a_{k-1}\}$ and $n - 1 \notin \sigma^{(n)}\{a_0, a_1, \ldots, a_{k-1}\}$.

\vspace{3mm}

Immediately from the proof of Proposition 2 follows

\vspace{3mm}

\textbf{Corollary} \z  \textsl{For any internal simplex $\sigma^{(n)}\{a_0, a_1, \ldots, a_{k-1}\}$  semirings $\mathcal{DN}^{\,1}_m$ are commutative and all their elements are $a_m$--nilpotent, where $m = 0, \ldots, k-1$.}

\vspace{3mm}

\textbf{Proposition} \z  \textsl{Let $\overline{a_m}$, where $m = 0, \ldots, k-1$, be a vertex of  internal simplex $\sigma^{(n)}\{a_0, a_1, \ldots, a_{k-1}\}$.  Then the set $\mathcal{DN}^{\,2}_m = \mathcal{DN}^{\,1}_m\cup \mathcal{L}^{n-2}_{a_m}\left(\sigma^{(n)}\{a_0, a_1, \ldots, a_{k-1}\}\right)$, where $m = 0, \ldots, k-1$, is a  subsemiring of $\sigma^{(n)}\{a_0, a_1, \ldots, a_{k-1}\}$.}

\emph{Proof.} The elements of $\mathcal{DN}^{\,2}_m$ are: $\overline{a_m}$, $a_i(a_m)_{n-1}$, where $i = 0, \ldots, m-1$, $(a_m)_{n-1}a_j$, where $j = m+1, \ldots, k-1$, $a_pa_q(a_m)_{n-2}$, where $p, q = 0, \ldots, m-1$, $p \leq q$,  $(a_m)_{n-2}a_ra_s$, where $r, s = m+1, \ldots, k-1$, $r \leq s$, and $a_p(a_m)_{n-2}a_s$, where $p = 0, \ldots, m-1$, $s = m+1, \ldots, k-1$.

Since $\mathcal{DN}^{\,1}_m$ is closed under the addition in order to prove the same for $\mathcal{DN}^{\,2}_m$, we consider:
$$a_i(a_m)_{n-1}\! + a_pa_q(a_m)_{n-2} =\! \!\left\{ \begin{array}{ll} a_pa_q(a_m)_{n-2} & \! \mbox{if}\; i \leq p\\
a_ia_q(a_m)_{n-2} & \! \mbox{if}\; i > p \end{array}\! \right.\!, a_i(a_m)_{n-1}\! + (a_m)_{n-2}a_ra_s\! =\! (a_m)_{n-2}a_ra_s,$$
$$(a_m)_{n-1}a_j\! + (a_m)_{n-2}a_ra_s =\! \!\left\{ \begin{array}{ll} (a_m)_{n-2}a_ra_s & \! \mbox{if}\; j \leq s\\
(a_m)_{n-2}a_ra_j & \! \mbox{if}\; j > s \end{array}\! \right.\!, (a_m)_{n-1}a_j\! + a_pa_q(a_m)_{n-2}\! =\! (a_m)_{n-1}a_j,$$
$$a_i(a_m)_{n-1} + a_p(a_m)_{n-2}a_s = \left\{ \begin{array}{ll} a_p(a_m)_{n-2}a_s &  \mbox{if}\; i \leq p\\
a_i(a_m)_{n-2}a_s &  \mbox{if}\; i > p \end{array} \right.,$$
$$(a_m)_{n-1}a_j + a_p(a_m)_{n-2}a_s = \left\{ \begin{array}{ll} a_p(a_m)_{n-2}a_s &  \mbox{if}\; j \leq s\\
a_p(a_m)_{n-2}a_j &  \mbox{if}\; j > s \end{array} \right.,$$
$$a_pa_q(a_m)_{n-2} + (a_m)_{n-2}a_ra_s = (a_m)_{n-2}a_ra_s,$$
$$a_{p_0}a_{q_0}(a_m)_{n-2} + a_p(a_m)_{n-2}a_s = \left\{ \begin{array}{ll} a_p(a_m)_{n-2}a_s &  \mbox{if}\; p_0 \leq p\\
a_{p_0}(a_m)_{n-2}a_s &  \mbox{if}\; p_0 > p \end{array} \right.,$$
$$(a_m)_{n-2}a_{r_0}a_{s_0} + a_p(a_m)_{n-2}a_s = \left\{ \begin{array}{ll} a(a_m)_{n-2}a_{r_0}a_s &  \mbox{if}\; s_0 \leq s\\
(a_m)_{n-2}a_{r_0}a_{s_0} &  \mbox{if}\; s_0 > s \end{array} \right.,$$
$$ a_p(a_m)_{n-2}a_s + a_{p_0}(a_m)_{n-2}a_{s_0} = \left\{ \begin{array}{ll} a_p(a_m)_{n-2}a_s &  \mbox{if}\; p \leq p_0, s \leq s_{0}\\
a_p(a_m)_{n-2}a_{s_0} &  \mbox{if}\; p \leq p_0, s > s_{0}\\
a_{p_0}(a_m)_{n-2}a_s &  \mbox{if}\; p > p_0, s \leq s_{0}\\
a_{p_0}(a_m)_{n-2}a_{s_0} &  \mbox{if}\; p > p_0, s > s_{0}\\
\end{array} \right.,$$
$$ \overline{a_m} + a_pa_q(a_m)_{n-2} = \overline{a_m},\; \overline{a_m} + (a_m)_{n-2}a_ra_s = (a_m)_{n-2}a_ra_s,\; \overline{a_m} + a_p(a_m)_{n-2}a_s = (a_m)_{n-1}a_s,$$
 where $i, p, q, p_0, q_0 = 0, 1, \ldots, m-1$, $p \leq q$, $p_0 < q_0$, $j, r, s,r_0, s_0 = {m+1}, \ldots, {k-1}$, $r \leq s$, $r_0 < s_0$. So, we prove that $\mathcal{DN}^{\,2}_m$ is closed under the addition.

Now we consider three cases, where, for the indices,  the upper restrictions are fulfilled.


\emph{Case 1.} Let $a_m = 1$. We shall show that all endomorphisms of $\mathcal{DN}^{\,2}_1$ are $1$--nilpotent with the only exception when $a_{k-1} = n-2$. When $a_{k-1} < n-2$, since $1$ is the least image of any endomorphism, there are only a few equalities: $ 1_{n-2}a_ra_s\cdot 1_{n-2}a_{r_0}a_{s_0} = \overline{1}$,
$$ 1_{n-1}a_j \cdot 1_{n-2}a_ra_s = 1_{n-2}a_ra_s \cdot 1_{n-1}a_j = \overline{1},\;  \overline{1}\cdot 1_{n-2}a_ra_s = 1_{n-2}a_ra_s \cdot \overline{1} = \overline{1}.$$
Hence, it follows that $\mathcal{DN}^{\,2}_1$ is a commutative semiring with trivial multiplication.

If $a_{k-1} = n-2$ it is easy to see that endomorphism $1_{n-2}(n-2)_2$ is the unique idempotent of $\mathcal{DN}^{\,2}_1$ (see [2]). Now we find
$1_{n-2}(n-2)_2 \cdot 1_{n-2}a_ra_s = 1_{n-2}(a_r)_2$,  $1_{n-1}(n-2) \cdot 1_{n-2}a_ra_s = 1_{n-1}a_r$, $1_{n-2}a_ra_s \cdot 1_{n-1}(n-2) = \overline{1}$.
Hence, $\mathcal{DN}^{\,2}_1$ is a semiring.


\emph{Case 2.} Let $a_m = n-2$. We shall show that all the endomorphisms of $\mathcal{DN}^{\,2}_{n-2}$ are $1$--nilpotent with the only exception when $a_0 = 1$.
 When $a_0 > 1$ we find: $$ a_pa_q(n-2)_{n-2}\cdot a_{p_0}a_{q_0}(n-2)_{n-2} = \overline{n-2},$$
 $$a_i(n-2)_{n-1}\cdot  a_pa_q(n-2)_{n-2} =  a_pa_q(n-2)_{n-2}\cdot a_i(n-2)_{n-1} = \overline{n-2},$$
 $$\overline{n-2}\cdot a_pa_q(n-2)_{n-2} = a_pa_q(n-2)_{n-2}\cdot \overline{n-2} = \overline{n-2}.$$

 If $a_0 = 1$ the only idempotent is $1_2(n-2)_{n-2}$ and we find:
$$1_2(n-2)_{n-2}\cdot a_pa_q(n-2)_{n-2} = (a_q)_2(n-2)_{n-2},\;$$
$$1(n-2)_{n-1}\cdot  a_pa_q(n-2)_{n-2} = a_q(n-2)_{n-1},\; a_pa_q(n-2)_{n-2}\cdot 1(n-2)_{n-1} = \overline{n-2},$$
Hence, $\mathcal{DN}^{\,2}_{n-2}$ is a semiring.

\emph{Case 3.}  Let $1 < a_0$ and $a_{k-1} < n-2$. We find the following trivial equalities, which are grouped by duality:
$$a_pa_q(a_m)_{n-2}\cdot a_{p_0}a_{q_0}(a_m)_{n-2} = \overline{a_m},\; (a_m)_{n-2}a_ra_s\cdot (a_m)_{n-2}a_{r_0}a_{s_0} = \overline{a_m},$$
$$a_pa_q(a_m)_{n-2}\cdot a_{p_0}(a_m)_{n-2}a_{s_0} = a_{p_0}(a_m)_{n-2}a_{s_0} \cdot a_pa_q(a_m)_{n-2} = \overline{a_m},$$
$$(a_m)_{n-2}a_ra_s\cdot a_{p_0}(a_m)_{n-2}a_{s_0} = a_{p_0}(a_m)_{n-2}a_{s_0} \cdot (a_m)_{n-2}a_ra_s = \overline{a_m},$$
$$a_pa_q(a_m)_{n-2}\cdot (a_m)_{n-2}a_ra_{s} = (a_m)_{n-2}a_ra_{s} \cdot a_pa_q(a_m)_{n-2} = \overline{a_m},$$
$$a_i(a_m)_{n-1}\cdot a_pa_q(a_m)_{n-2} = a_pa_q(a_m)_{n-2}\cdot a_i(a_m)_{n-1} = \overline{a_m},$$
$$(a_m)_{n-1}a_j\cdot a_pa_q(a_m)_{n-2}  = a_pa_q(a_m)_{n-2}\cdot (a_m)_{n-1}a_j = \overline{a_m},$$
$$a_i(a_m)_{n-1}\cdot a_p(a_m)_{n-2}a_s = a_p(a_m)_{n-2}a_s\cdot a_i(a_m)_{n-1} = \overline{a_m},$$
$$(a_m)_{n-1}a_j\cdot a_p(a_m)_{n-2}a_s  = a_p(a_m)_{n-2}a_s\cdot (a_m)_{n-1}a_j = \overline{a_m},$$
$$a_i(a_m)_{n-1}\cdot (a_m)_{n-2}a_ra_s = (a_m)_{n-2}a_ra_s\cdot a_i(a_m)_{n-1} = \overline{a_m},$$
$$(a_m)_{n-1}a_j\cdot (a_m)_{n-2}a_ra_s  = a(a_m)_{n-2}a_ra_s\cdot (a_m)_{n-1}a_j = \overline{a_m},$$
$$\overline{a_m}\cdot a_pa_q(a_m)_{n-2} = a_pa_q(a_m)_{n-2}\cdot \overline{a_m} = \overline{a_m},$$
$$\overline{a_m}\cdot a_p(a_m)_{n-2}a_s = a_p(a_m)_{n-2}a_s\cdot \overline{a_m} = \overline{a_m},$$
$$\overline{a_m}\cdot (a_m)_{n-2}a_ra_s = (a_m)_{n-2}a_ra_s\cdot \overline{a_m} = \overline{a_m}.$$

\emph{Case 4.}  Let $a_0 = 1$ and $a_{k-1} < n-2$. Then  $\;1_2(a_m)_{n-2}$ is the only idempotent in $\mathcal{DN}^{\,2}_m$. Additionally to the equalities of the previous case we find:
$$1_2(a_m)_{n-2}\cdot a_pa_q(n-2)_{n-2} = (a_q)_2(a_m)_{n-2}, \; 1(a_m)_{n-1}\cdot a_pa_q(a_m)_{n-2} = a_q(a_m)_{n-1}.$$


\emph{Case 5.}  Let $1 < a_0$ and $a_{k-1} = n-2$. Now the only idempotent endomorphism in $\mathcal{DN}^{\,2}_m$ is $(a_m)_{n-2}(n-2)_2$. We additionally find the following equalities:
$$(a_m)_{n-2}(n-2)_2 \cdot (a_m)_{n-2}a_ra_s = (a_m)_{n-2}(a_r)_2,\; (a_m)_{n-1}(n-2)\cdot (a_m)_{n-2}a_ra_s  = (a_m)_{n-1}a_r.$$


\emph{Case 6.}  Let $a_0 = 1$ and $a_{k-1} = n-2$. Now, in $\mathcal{DN}^{\,2}_m$, there are two idempotents: $\;1_2(a_m)_{n-2}$ and $(a_m)_{n-2}(n-2)_2$. Here the equalities from cases 4 and 5 are valid and also all the equalities from case 3, under the respective restrictions for the indices, are fulfilled.

Hence, $\mathcal{DN}^{\,2}_m$ is a semiring. \hfill $\Box$

\vspace{2mm}

\textbf{Theorem} \z  \textsl{Let $\sigma^{(n)}_k(A) = \sigma^{(n)}\{a_0, a_1, \ldots, a_{k-1}\}$ be a simplex. }

\textsl{a. For the least vertex $\overline{a_0}$ it follows $\mathcal{DN}^{\,n-a_0-1}_0 = \sigma^{(n)}_k(A)\cap {\mathcal{E}}^{(a_0)}_{\mathcal{C}_n}$.
}

\textsl{b. For the biggest vertex $\overline{a_{k-1}}$ it follows $\mathcal{DN}^{\,a_{k-1}}_{k-1} = \sigma^{(n)}_k(A)\cap {\mathcal{E}}^{(a_{k-1})}_{\mathcal{C}_n}$.
}

\emph{Proof.} a.
 Since $\overline{a_0}$ is the least vertex of the simplex, it follows that  layer
$\mathcal{L}^{a_0 + 1}_{a_0}\left(\sigma^{(n)}\{a_0, a_1, \ldots, a_{k-1}\}\right)$ consists of endomorphisms
${\alpha = (a_0)_{a_0 + 1}(a_1)_{p_1} \ldots (a_{k-1})_{p_{k-1}}}$, where $a_0 + 1 + p_1 + \cdots + p_{k-1} = n$, i.e. $\alpha(0) = a_0$, $\ldots$, $\alpha(a_0) = a_0$. All  the layers
$\mathcal{L}^{\ell}_{a_0}\left(\sigma^{(n)}\{a_0, a_1, \ldots, a_{k-1}\}\right)$, where $\ell \geq a_0 + 1$, consist of endomorphisms having $a_0$ as a fixed point. So, $\mathcal{DN}^{\,n-a_0-1}_0 \subseteq \sigma^{(n)}_k(A)\cap {\mathcal{E}}^{(a_0)}_{\mathcal{C}_n}$.

Conversely, let $\alpha \in \sigma^{(n)}_k(A)\cap {\mathcal{E}}^{(a_0)}_{\mathcal{C}_n}$. Then $\alpha(a_0) = a_0$. Since $\overline{a_0}$ is the least vertex of the simplex, we have $\alpha(0) = \ldots = \alpha(a_0 - 1) = a_0$, that is $\alpha \in \mathcal{L}^{\ell}_{a_0}\left(\sigma^{(n)}\{a_0, a_1, \ldots, a_{k-1}\}\right)$, where $\ell \geq a_0 + 1$. Hence, $\mathcal{DN}^{\,n-a_0-1}_0 = \sigma^{(n)}_k(A)\cap {\mathcal{E}}^{(a_0)}_{\mathcal{C}_n}$.

b.  Since $\overline{a_{k-1}}$ is the biggest vertex of the simplex, it follows that  layer
$\mathcal{L}^{n - a_{k-1}}_{a_{k-1}}\left(\sigma^{(n)}\{a_0, a_1, \ldots, a_{k-1}\}\right)$ consists of endomorphisms ${\alpha = (a_0)_{p_0} \ldots (a_{k-2})_{p_{k-2}}}(a_{k-1})_{n- a_{k-1}}$, where $p_0 +  \cdots + p_{k-2} + n - a_{k-1} = n$. So, $p_0 +  \cdots + p_{k-2} = a_{k-1}$ implies that the images of $0$, $\ldots$, $a_{k-1} - 1$ are not equal to $a_{k-1}$, but $\alpha(a_{k-1}) = a_{k-1}$.
For all the endomorphisms of  layers
$\mathcal{L}^{\ell}_{a_{k-1}}\left(\sigma^{(n)}\{a_0, a_1, \ldots, a_{k-1}\}\right)$, where $\ell \geq n - a_{k-1}$, we have $p_0 +  \cdots + p_{k-2} = a_{k-1}$. Hence, the elements of these layers have $a_{k-1}$ as a fixed point and $\mathcal{DN}^{\,a_{k-1}}_0 \subseteq \sigma^{(n)}_k(A)\cap {\mathcal{E}}^{(a_{k-1})}_{\mathcal{C}_n}$.

Conversely, let $\alpha \in \sigma^{(n)}_k(A)\cap {\mathcal{E}}^{(a_{k-1})}_{\mathcal{C}_n}$. Then $\alpha(a_{k-1}) = a_{k-1}$. Since $\overline{a_{k-1}}$ is the biggest vertex of the simplex, we have $\alpha(a_{k-1} + 1) = \ldots = \alpha(n - 1) = a_{k-1}$. Thus,  $\alpha \in \mathcal{L}^{\ell}_{a_{k-1}}\left(\sigma^{(n)}\{a_0, a_1, \ldots, a_{k-1}\}\right)$, where $\ell \geq n - a_{k-1}$. Hence, $\mathcal{DN}^{\,a_{k-1}}_{k-1} = \sigma^{(n)}_k(A)\cap {\mathcal{E}}^{(a_{k-1})}_{\mathcal{C}_n}$. \hfill $\Box$

\vspace{2mm}

\emph{{Remark}} \z What is the least $\ell$, such that the discrete neighborhood $\mathcal{DN}^{\,\ell}_{m}$ of the vertex $\overline{a_m}$ of simplex $\sigma^{(n)}\{a_0, a_1, \ldots, a_{k-1}\}$, where $m \neq 0$ and $m \neq k-1$, is a semiring? Since $1_32_{n-2}$ is an $1$--nilpotent element of any simplex
$\sigma^{(n)}\{1, 2, a_2 \ldots, a_{k-1}\}$, it follows that $\ell = 2$.

\vspace{1mm}

From the last theorem it follows that all the $a_0$--nilpotent elements of the simplex $\sigma^{(n)}\{a_0, a_1, \ldots, a_{k-1}\}$ are from  semiring $\mathcal{DN}^{\,n-a_0-1}_0$. But there are elements of $\mathcal{DN}^{\,n-a_0-1}_0$ which are not $a_0$--nilpotent. For instance,  endomorphism $\alpha \in \mathcal{L}^{a_0 + 1}_{a_0}\left(\sigma^{(n)}\{a_0, a_1, \ldots, a_{k-1}\}\right)$ such that $\alpha(i) = a_m$, where $m = 1, \ldots, k-1$, for any $i > a_0$ is an idempotent. In order to separate  $a_0$--nilpotent elements from all the other elements of $\mathcal{DN}^{\,n-a_0-1}_0$,  we consider the following

\vspace{2mm}

\textbf{Proposition} \z  \textsl{The endomorphism $\alpha \in \mathcal{DN}^{\,n-a_0-1}_0$  is $a_0$--nilpotent if}
$$\alpha(0) = \cdots = \alpha(a_0) = \cdots = \alpha(a_1) = a_0, \; \alpha(i) < i, \; \mbox{\textsl{for}}\; a_1 < i \leq n - 1.$$

\emph{Proof.} Let us suppose that for some $i \geq a_0 + 1$ follows $\alpha(i) \geq i$. Then  $\alpha^m(i) \geq i \geq a_0 + 1$ for any natural $m$, which contradicts that  $\alpha$ is $a_0$--nilpotent endomorphism. Hence, $\alpha(i) < i$ for $i \geq a_0 + 1$. In particular $\alpha(a_s) < a_s$, for any $s = 1, \ldots, k$ and then $\alpha(a_1) = a_0$. \hfill $\Box$

\vspace{1mm}

From the last proposition it immediately follows that there are not $a_0$--nilpotent endomorphisms in  layer $\mathcal{L}^{a_0 + 1}_{a_0}\left(\sigma^{(n)}\{a_0, a_1, \ldots, a_{k-1}\}\right)$. Since $\mathcal{DN}^{\,n-a_0-1}_0$ is a proper subsemiring of the simplex, it follows that in  layer
 $\mathcal{L}^{a_0 + 1}_{a_0}\left(\sigma^{(n)}\{a_0, a_1, \ldots, a_{k-1}\}\right)$ there are not any $a_m$--nilpotent elements. So, the elements of this layer are idempotents or roots of idempotents. But if $\alpha = (a_0)_{a_0 + 1}(a_1)_{p_1} \ldots (a_{k-1})_{p_{k-1}} \in \mathcal{L}^{a_0 + 1}_{a_0}\left(\sigma^{(n)}\{a_0, a_1, \ldots, a_{k-1}\}\right)$, it follows by induction that $\alpha^m = (a_0)_{a_0 + 1}(a_1)_{q_1} \ldots (a_{k-1})_{q_{k-1}} \in \mathcal{L}^{a_0 + 1}_{a_0}\left(\sigma^{(n)}\{a_0, a_1, \ldots, a_{k-1}\}\right)$ for any natural $m$. So, all the elements of this layer are idempotents or roots of idempotents of the same layer. Obviously, the layer is closed under the addition. So, we prove

\vspace{3mm}

\textbf{Proposition} \z  \textsl{For any simplex $\sigma^{(n)}\{a_0, a_1, \ldots, a_{k-1}\}$  layer\\ $\mathcal{L}^{a_0 + 1}_{a_0}\left(\sigma^{(n)}\{a_0, a_1, \ldots, a_{k-1}\}\right)$ is a subsemiring of the simplex.}

\vspace{5mm}

 \noindent{\bf \large 3. \hspace{0.5mm} Strings}

\vspace{3mm}

Let us denote the  elements of  semiring $\mathcal{STR}^{(n)}\{a,b\}$ by $a_kb_{n-k}$, where $k = 0,\ldots,n$ is the number of the elements of $\mathcal{C}_n$ with an image equal to $a$, i.e.
$$a_kb_{n-k} = \wr\, \underbrace{a, \ldots, a}_{k}, \underbrace{b, \ldots, b}_{n-k}\,\wr.$$

In particular, we denote $a_nb_0 = \overline{a}$ and $a_0b_n = \overline{b}$.

\vspace{3mm}

Let us consider the following subset of $\mathcal{STR}^{(n)}\{a,b\}$:
$$N^{[a]}\left(\mathcal{STR}^{(n)}\{a,b\}\right) = \{\overline{a}, \ldots, a_{b+1}b_{n-b-1}\}.$$
For any endomorphism $\alpha \in N^{[a]}\left(\mathcal{STR}^{(n)}\{a,b\}\right)$, that is $\alpha = a_{\ell\,} b_{n-\ell}$, where $b + 1 \leq \ell \leq n$, we have $\alpha(b) = a$. Hence $\alpha, \beta \in N^{[a]}\left(\mathcal{STR}^{(n)}\{a,b\}\right)$ implies $\alpha\cdot\beta = \overline{a}$. Since $\alpha^2 = \overline{a}$ for any $\alpha \in N^{[a]}\left(\mathcal{STR}^{(n)}\{a,b\}\right)$, it follows, see [3], that $$N^{[a]}\left(\mathcal{STR}^{(n)}\{a,b\}\right) = N^{[a]}_n\cap \mathcal{STR}^{(n)}\{a,b\}.$$

From Theorem 3.3 of [3], see section 1, follows

\vspace{3mm}

\textbf{Proposition} \z  \textsl{The set $N^{[a]}\left(\mathcal{STR}^{(n)}\{a,b\}\right)$ is a subsemiring of $\mathcal{STR}^{(n)}\{a,b\}$ consisting of  all $a$ -- nilpotent elements of this string.}

The order of this semiring is $n-b$.

\vspace{3mm}

The next subset of $\mathcal{STR}^{(n)}\{a,b\}$ is:
$$Id\left(\mathcal{STR}^{(n)}\{a,b\}\right) = \{a_bb_{n-b}, \ldots, a_{a+1}b_{n-a-1}\}.$$
For any endomorphism of $Id\left(\mathcal{STR}^{(n)}\{a,b\}\right)$  elements $a$ and $b$ are fixed points. From Corollary 3 of [2] we find that all the elements of $Id\left(\mathcal{STR}^{(n)}\{a,b\}\right)$ are idempotents and from Theorem 9 of [2], see section 1, it  follows

\vspace{3mm}

\textbf{Proposition} \z  \textsl{The set $Id\left(\mathcal{STR}^{(n)}\{a,b\}\right)$ is a subsemiring of $\mathcal{STR}^{(n)}\{a,b\}$ consisting of all idempotent elements of this string different from $\overline{a}$ and $\overline{b}$.
The order of this semiring is $b-a$.}

\vspace{3mm}

The last considered subset of $\mathcal{STR}^{(n)}\{a,b\}$ is:
$$N^{[b]}\left(\mathcal{STR}^{(n)}\{a,b\}\right) = \{a_ab_{n-a}, \ldots, \overline{b}\}.$$
For any endomorphism $\alpha = a_{\ell\,} b_{n-\ell}$, where $0 \leq \ell \leq a$ it follows $\alpha(a) = b$. Hence $\alpha, \beta \in N^{[b]}\left(\mathcal{STR}^{(n)}\{a,b\}\right)$ implies $\alpha\cdot\beta = \overline{b}$. Since $\alpha^2 = \overline{b}$ for any $\alpha \in N^{[b]}\left(\mathcal{STR}^{(n)}\{a,b\}\right)$ it follows, see [2], that $$N^{[b]}\left(\mathcal{STR}^{(n)}\{a,b\}\right) = N^{[b]}_n\cap \mathcal{STR}^{(n)}\{a,b\}.$$

From Theorem 3.3 of [3], see section 1, we have

\vspace{3mm}

\textbf{Proposition} \z  \textsl{The set $N^{[b]}\left(\mathcal{STR}^{(n)}\{a,b\}\right)$ is a subsemiring of $\mathcal{STR}^{(n)}\{a,b\}$ consisting of all $b$ -- nilpotent elements of this string.}

The order of this semiring is $a + 1$.

\vspace{3mm}

\textbf{Proposition} \z  \textsl{Let $a, b \in \mathcal{C}_n$, $a < b$ and $a_kb_{n-k} \in \mathcal{STR}^{(n)}\{a,b\}$, where $k = 0,\ldots,n$. Then}
$$\begin{array}{ll}
a_kb_{n-k} \cdot \alpha = \overline{a}, & \mbox{\textsl{if}} \;\; \alpha \in N^{[a]}\left(\mathcal{STR}^{(n)}\{a,b\}\right)\\
a_kb_{n-k} \cdot \alpha = a_kb_{n-k}, & \mbox{\textsl{if}} \;\; \alpha \in Id\left(\mathcal{STR}^{(n)}\{a,b\}\right)\\
a_kb_{n-k} \cdot \alpha = \overline{b}, & \mbox{\textsl{if}} \;\; \alpha \in N^{[b]}\left(\mathcal{STR}^{(n)}\{a,b\}\right)
\end{array}.$$

\emph{Proof.} For any $i \in \mathcal{C}_n$ and $\alpha \in N^{[a]}\left(\mathcal{STR}^{(n)}\{a,b\}\right)$ it follows
$$(a_kb_{n-k} \cdot \alpha)(i) = \alpha(a_kb_{n-k}(i)) = \left\{ \begin{array}{ll} \alpha(a),& \mbox{if}\;\; 0 \leq i \leq k\\ \alpha(b), & \mbox{if}\; \; k+1 \leq i \leq n - 1 \end{array} =\: a \right.$$
which means that $a_kb_{n-k} \cdot \alpha = \overline{a}$.

 For any $i \in \mathcal{C}_n$ and $\alpha \in Id \left(\mathcal{STR}^{(n)}\{a,b\}\right)$ it follows
$$(a_kb_{n-k} \cdot \alpha)(i) = \alpha(a_kb_{n-k}(i)) = \left\{ \begin{array}{ll} \alpha(a),& \mbox{if}\;\; 0 \leq i \leq k\\ \alpha(b), & \mbox{if}\; \; k+1 \leq i \leq n - 1 \end{array} = \right.$$ $$\hphantom{aaaaaaaaaaaaaaaaaaaaaaaaaaaaaaaaa} \left\{ \begin{array}{ll} a,& \mbox{if}\;\; 0 \leq i \leq k\\ b, & \mbox{if}\; \; k+1 \leq i \leq n - 1 \end{array} =  \: a_kb_{n-k}(i) \right.$$
which means that $a_kb_{n-k} \cdot \alpha = a_kb_{n-k}$.

For any $i \in \mathcal{C}_n$ and $\alpha \in N^{[b]}\left(\mathcal{STR}^{(n)}\{a,b\}\right)$ it follows
$$(a_kb_{n-k} \cdot \alpha)(i) = \alpha(a_kb_{n-k}(i)) = \left\{ \begin{array}{ll} \alpha(a),& \mbox{if}\;\; 0 \leq i \leq k\\ \alpha(b), & \mbox{if}\; \; k+1 \leq i \leq n - 1 \end{array} =\: b \right.$$
which means that $a_kb_{n-k} \cdot \alpha = \overline{b}$. \hfill $\Box$

\vspace{3mm}

Immediately follows

\vspace{3mm}

\textbf{Corollary} \z  \textsl{The idempotent endomorphisms of  semiring $\mathcal{STR}^{(n)}\{a,b\}$, different from $\overline{a}$ and $\overline{b}$, are right identities.}

\vspace{3mm}

\textbf{Corollary} \z  \textsl{Any two different strings are nonisomorphic semirings.}

\vspace{3mm}

Using the fact that the strings are faces of any $k$ -- simplex for arbitrary $k \geq 2$,  the last corollary implies

\vspace{3mm}

\textbf{Corollary} \z  \textsl{Any two different $k$ -- simplices are nonisomorphic semirings.}

\vspace{3mm}

\emph{{Remark}} \z a. From Propopsition 1.4 we actually observe that the multiplicative structure of arbitrary string $\mathcal{STR}^{(n)}\{a,b\}$ is very clear: first, we find $n - b$ endomorphisms (all the $a$ -- nilpotent elements) which are square roots of $\overline{a}$ or $\overline{a}$, then $b - a$ idempotents, which are right identities and, in the end, $a + 1$ elements (all the $b$ -- nilpotent elements) which are square roots of $\overline{b}$ or $\overline{b}$.

b. The union of semirings $N^{[a]}\left(\mathcal{STR}^{(n)}\{a,b\}\right)$ and $Id\left(\mathcal{STR}^{(n)}\{a,b\}\right)$ is also a semi\-ring because
$$N^{[a]}\left(\mathcal{STR}^{(n)}\{a,b\}\right)\cup Id\left(\mathcal{STR}^{(n)}\{a,b\}\right) = \mathcal{STR}^{(n)}\{a,b\}\cap {\mathcal{E}}^{(a)}_{\mathcal{C}_n}.$$

Similarly,
$$N^{[b]}\left(\mathcal{STR}^{(n)}\{a,b\}\right)\cup Id\left(\mathcal{STR}^{(n)}\{a,b\}\right) = \mathcal{STR}^{(n)}\{a,b\}\cap {\mathcal{E}}^{(b)}_{\mathcal{C}_n}$$
is a subsemiring of $\mathcal{STR}^{(n)}\{a,b\}$.

\vspace{3mm}

Two strings $\mathcal{STR}^{(n)}\{a,b\}$ and $\mathcal{STR}^{(n)}\{x,y\}$ are called {\emph{consecutive}}  if they have a common vertex. So,  strings $\mathcal{STR}^{(n)}\{a,b\}$ and $\mathcal{STR}^{(n)}\{b,c\}$, $\mathcal{STR}^{(n)}\{a,b\}$ and $\mathcal{STR}^{(n)}\{a,c\}$, $\mathcal{STR}^{(n)}\{a,c\}$ and $\mathcal{STR}^{(n)}\{b,c\}$ (when $a < b < c$) are the three possibilities of the pairs of consecutive strings.

Let $a_kb_{n-k} \in \mathcal{STR}^{(n)}\{a,b\}$, where $k = 0,\ldots,n$, and $b_\ell c_{n-\ell} \in \mathcal{STR}^{(n)}\{b,c\}$, where $\ell = 0,\ldots,n$. Since $a_kb_{n-k} < b_\ell c_{n-\ell}$, then $a_kb_{n-k} + b_\ell c_{n-\ell} = b_\ell c_{n-\ell}$. By similar arguments, for any $a_m c_{n-m} \in \mathcal{STR}^{(n)}\{a,c\}$, we can construct $a_m c_{n-m} + b_\ell c_{n-\ell} = b_rc_{n-r}$, where $r = \min\{\ell,m\}$. But when we add  endomorphisms $a_k b_{n-k}$ and $a_mc_{n-m}$, where $k < m$, the sum is $a_kb_{m-k}c_{n-m}$, so, the set of these three strings is not closed under the addition.

In the next proposition we examine the product of endomorphisms of two (not necessarily consecutive) strings.

\vspace{3mm}

\textbf{Proposition} \z \textsl{Let $a_kb_{n-k} \in \mathcal{STR}^{(n)}\{a,b\}$, where $k = 0,\ldots,n$, and $x_\ell\, y_{n-\ell} \in \mathcal{STR}^{(n)}\{x,y\}$, where $\ell = 0,\ldots,n$. Then}
$$\begin{array}{ll}
a_kb_{n-k} \cdot x_\ell\, y_{n-\ell} = \overline{x}, & \mbox{\textsl{if}} \;\; b+1 \leq \ell \leq n\\
a_kb_{n-k} \cdot x_\ell\, y_{n-\ell} = x_ky_{n-k}, & \mbox{\textsl{if}} \;\; a + 1 \leq \ell \leq b\\
a_kb_{n-k} \cdot x_\ell\, y_{n-\ell} = \overline{y}, & \mbox{\textsl{if}} \;\; 0 \leq \ell \leq a
\end{array}.$$

\emph{Proof.} For any $i \in \mathcal{C}_n$ it follows
$$\left(a_kb_{n-k} \cdot x_\ell\, y_{n-\ell}\right)(i) = x_\ell\, y_{n-\ell}\left(a_kb_{n-k}\right)(i) = \left\{ \begin{array}{ll} x_\ell\, y_{n-\ell}(a),& \mbox{if}\;\; 0 \leq i \leq k\\ x_\ell\, y_{n-\ell}(b), & \mbox{if}\; \; k+1 \leq i \leq n - 1 \end{array}. \right.$$

If $b+1 \leq \ell \leq n$, we have $x_\ell\, y_{n-\ell}(b) = x$ and then $x_\ell\, y_{n-\ell}(a) = x$. So, for any $i \in \mathcal{C}_n$ we find that $\left(a_kb_{n-k} \cdot x_\ell\, y_{n-\ell}\right)(i) = x$ and hence $a_kb_{n-k} \cdot x_\ell\, y_{n-\ell} = \overline{x}$.

If $a + 1 \leq \ell \leq b$, it follows $x_\ell\, y_{n-\ell}(a) = x$ and $x_\ell\, y_{n-\ell}(b) = y$. Thus, we obtain
$$\left(a_kb_{n-k} \cdot x_\ell\, y_{n-\ell}\right)(i) =  \left\{ \begin{array}{ll} x,& \mbox{if}\;\; 0 \leq i \leq k\\ y, & \mbox{if}\; \; k+1 \leq i \leq n - 1 \end{array} \right. = x_ky_{n-k}(i).$$
Hence, $a_kb_{n-k} \cdot x_\ell\, y_{n-\ell} = x_ky_{n-k}$.
If $ 0 \leq \ell \leq a$, we have $x_\ell\, y_{n-\ell}(a) = y$ and then $x_\ell\, y_{n-\ell}(b) = y$. So, for any $i \in \mathcal{C}_n$ we obtain $\left(a_kb_{n-k} \cdot x_\ell\, y_{n-\ell}\right)(i) = y$ and hence $a_kb_{n-k} \cdot x_\ell\, y_{n-\ell} = \overline{y}$. \hfill $\Box$

\vspace{3mm}

Immediately it follows

\vspace{3mm}

\textbf{Corollary} \z  \textsl{For any $a, b, c \in \mathcal{C}_n$, $a < b < c$, the set consisting of all the elements of the consecutive strings $\mathcal{STR}^{(n)}\{a,b\}$ and $\mathcal{STR}^{(n)}\{b,c\}$ is a semiring.}

\vspace{3mm}

A subsemiring $S$ of  endomorphism semiring $\widehat{\mathcal{E}}_{\mathcal{C}_n}$ is called a {\emph{trivial semiring}} if  for any two elements $\alpha, \beta \in S$ it follows $\alpha\cdot \beta = \iota$, where $\iota$ is a fixed element of $S$. If semiring $S$ is a trivial, then there exists a unique idempotent $\iota\in S$ such that the product of any two elements of $S$ is equal to $\iota$. If this idempotent is the biggest (least) element of the trivial semiring $S$, the $S$ is called  an \emph{upper (lower) trivial semiring}.

\vspace{3mm}

\emph{{Example}} \z  a) The semiring of $a$ -- nilpotent elements $N^{[a]}\left(\mathcal{STR}^{(n)}\{a,b\}\right)$ (using the proof of Proposition 1.4) is a lower trivial semiring.

b) The semiring of $b$ -- nilpotent elements $N^{[b]}\left(\mathcal{STR}^{(n)}\{a,b\}\right)$ (using the proof of Proposition 1.4) is an upper trivial semiring.

c) Let us consider the semiring from Corollary 2.2. Then from Proposition 2.1 follows that the union of  semirings $N^{[b]}\left(\mathcal{STR}^{(n)}\{a,b\}\right)$ and $N^{[b]}\left(\mathcal{STR}^{(n)}\{b,c\}\right)$ is a trivial semiring. Since $\overline{b}$ is the biggest element of $N^{[b]}\left(\mathcal{STR}^{(n)}\{a,b\}\right)$ and the least element of $N^{[b]}\left(\mathcal{STR}^{(n)}\{b,c\}\right)$ it follows that the considered trivial semiring is neither  upper trivial nor  lower trivial.

\vspace{3mm}
Now we shall construct some useful subsemirings of a given string, some of which are trivial semirings.
We consider the following subset of  string $\mathcal{STR}^{(n)}\{a,b\}$:
$$A_r = \{\overline{a}, a_{n-1}b, \ldots, a_rb_{n-r}\},$$
where $r = 1, \ldots, n$. Since $A_r$ is a chain for any $r$, it is closed under the addition. If $r \geq b + 1$, then $A_r \subseteq N^{[a]}\left(\mathcal{STR}^{(n)}\{a,b\}\right)$, so, $A_r$ is a lower trivial semiring.  If $r \leq a$, then $A_r \cap N^{[b]}\left(\mathcal{STR}^{(n)}\{a,b\}\right) \neq \varnothing$ which implies that $A_r$ is not closed under the multiplication, i.e. it is not a semiring.
From Remark 1.7 b. it follows that
$$A_{a+1} = N^{[a]}\left(\mathcal{STR}^{(n)}\{a,b\}\right)\cup Id\left(\mathcal{STR}^{(n)}\{a,b\}\right)
$$
is the biggest between  sets $A_r$, which is a semiring.
Since every element of semiring $Id\left(\mathcal{STR}^{(n)}\{a,b\}\right)$ is a right identity of  string $\mathcal{STR}^{(n)}\{a,b\}$, it follows that for every $r = a + 1, \ldots, n$ the set $A_r$ is a semiring.

Using the same idea, we consider the subset of $\mathcal{STR}^{(n)}\{a,b\}$:
$$B_s = \{\overline{b}, ab_{n-1}, \ldots, a_sb_{n-s}\},$$
where $s = 0, \ldots, n-1$. The set $B_s$ is a chain for any $s$, so, it is closed under the addition. If $s \leq a$, then $B_s \subseteq N^{[b]}\left(\mathcal{STR}^{(n)}\{a,b\}\right)$, so, $B_s$ is an upper trivial semiring. If $s \geq b+ 1$, then $B_s \cap N^{[a]}\left(\mathcal{STR}^{(n)}\{a,b\}\right) \neq \varnothing$ which means that $B_s$ is not closed under the multiplication, so, $B_s$ is not a semiring.
Also from Remark 1.7 b. it follows that
$$B_{b} = N^{[b]}\left(\mathcal{STR}^{(n)}\{a,b\}\right)\cup Id\left(\mathcal{STR}^{(n)}\{a,b\}\right)
$$
is the biggest from  sets $B_s$ which is a semiring.
 By the same way, considering  that every element of semiring $Id\left(\mathcal{STR}^{(n)}\{a,b\}\right)$ is a right identity of  string $\mathcal{STR}^{(n)}\{a,b\}$ it follows that for every $s = 0, \ldots, b$  set $B_s$ is a semiring.

\vspace{5mm}

\noindent{\large \bf 4. Triangles}

\vspace{3mm}

Let $a, b, c \in \mathcal{C}_n$, $a < b < c$, are fixed elements. The set of endomorphisms $\alpha$ such that
$$\alpha(0) = \cdots = \alpha(k-1) = a, \alpha(k) = \cdots = \alpha(k+\ell-1) = b, \alpha(k+\ell) = \cdots = \alpha(n-1) = c$$
or briefly $\alpha = a_kb_\ell c_{n-k-\ell}$, where $0 \leq k \leq n-1$, $0 \leq \ell \leq n-1$ and $0 \leq n - k -\ell \leq n-1$ is actually the triangle $\triangle^{(n)}\{a,b,c\}$. Obviously, the order of this semiring is $\displaystyle \binom{n+2}{2}$.

The strings $\mathcal{STR}^{(n)}\{a,b\}$, $\mathcal{STR}^{(n)}\{a,c\}$ and $\mathcal{STR}^{(n)}\{b,c\}$ are called strings of  $\triangle^{(n)}\{a,b,c\}$.

\vspace{3mm}

Let $R$ be a subsemiring of $\widehat{\mathcal{E}}_{\mathcal{C}_n}$ and $\alpha, \beta \in R$, $\alpha \neq \beta$. These endomorphisms are called {\emph{right-similar}} ({\emph{left-similar}}) if for any $\gamma \in R$ we have $\alpha\cdot \gamma = \beta\cdot \gamma$ ($\gamma\cdot \alpha = \gamma\cdot \beta$). We denote this by $\alpha \sim_r \beta$ ($\alpha \sim_\ell \beta$). In the next sections we shall answer  the question: Are there right-similar (left-similar) elements in $\triangle^{(n)}\{a,b,c\}$ ?

\vspace{3mm}

\emph{Example} \z The biggest side of the least tetrahedron $\mathcal{TETR}^{(4)}\{0,1,2,3\}$ is  triangle $\triangle^{(4)}\{1,2,3\}$. The elements of this semiring can be arranged as in the following scheme (fig.1):

\begin{figure}[h]\centering
  \includegraphics[width=67mm]{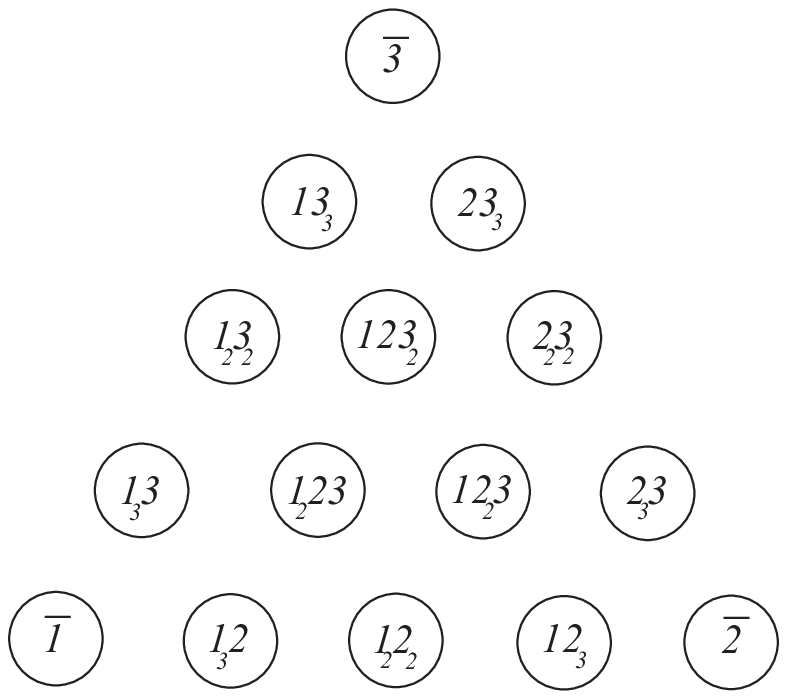}\\
\end{figure}
\vspace{-1mm}

\centerline{\small Figure 1.}
\vspace{4mm}

It is easy to see that the interior of this triangle is  set $Int = \{1_223, 12_23, 123_2\}$. Since $(123_2)^2 = 23_3 \notin Int$, it follows that the interior of $\triangle^{(4)}\{1,2,3\}$ is not a semiring. Since $1_223 = 1_22_2 + 1_33$, $12_23 = 12_3 + 1_33$ and $123_2 = 12_3 + 1_23_2$, it follows that every element of the interior of the triangle can be represented as a sum of an element of the least side of the triangle and an element of the middle side of the triangle.

In this triangle there are many left-similar endomorphisms: $12_3 \sim_\ell \overline{2}$, $13_3 \sim_\ell 23_3 \sim_\ell \overline{3}$, $123_2 \sim_\ell 2_23_2$, $12_23 \sim_\ell 2_33$. The endomorphism $1_223$ is a right-identity, so there are not right-similar endomorphisms in $\triangle^{(4)}\{1,2,3\}$.

\vspace{3mm}

\textbf{Proposition} \z \textsl{Any element of the interior of  $\triangle^{(n)}\{a,b,c\}$ can be uniquely represented as a sum of the elements of strings $\mathcal{STR}^{(n)}\{a,b\}$ and $\mathcal{STR}^{(n)}\{a,c\}$. }

\emph{Proof.} We easily calculate
$$a_{n-1}c + a_k b_{n-k} = a_k b_{n-k-1}c, \; \mbox{where} \; k = 0, \ldots, n-2.$$

By the same argument,  for any $j = 1, \ldots, n-1$ we find
$$a_{n-j}c_j + a_k b_{n-k} = a_k b_{n-k-j}c_j, \; \mbox{where} \; k = 0, \ldots, n-j-1.$$

So, we prove more: all the elements of the interior of $\triangle^{(n)}\{a,b,c\}$ and all the elements of the interior of $\mathcal{STR}^{(n)}\{b,c\}$ are sums of the elements of $\mathcal{STR}^{(n)}\{a,b\}$ and $\mathcal{STR}^{(n)}\{a,c\}$. From the construction we observe that endomorphisms $\overline{a}$  and $\overline{c}$ do not occur in these sums and every  representation of this type is unique. \hfill $\Box$

\vspace{3mm}

\textbf{Corollary} \z  \textsl{The boundary of an arbitrary triangle $\triangle^{(n)}\{a,b,c\}$ is a multiplicative semigroup but not a semiring.}

\emph{Proof.} Since the boundary of  $\triangle^{(n)}\{a,b,c\}$ is a union of strings $\mathcal{STR}^{(n)}\{a,b\}$, $\mathcal{STR}^{(n)}\{a,c\}$ and $\mathcal{STR}^{(n)}\{b,c\}$, from Proposition 8 it follows that this set is a multiplicative semigroup.  From the last proposition, it follows that the boundary of $\triangle^{(n)}\{a,b,c\}$ is not closed under the addition. \hfill $\Box$

\vspace{3mm}

\textbf{Corollary} \z  \textsl{The interior of an arbitrary triangle $\triangle^{(n)}\{a,b,c\}$, where $n \geq 4$,  is an additive semigroup but not a semiring.}

\emph{Proof.} From the last proposition and the fact that $\mathcal{STR}^{(n)}\{a,b\}$ and $\mathcal{STR}^{(n)}\{a,c\}$ are semirings, it follows that the interior of the triangle is closed under the addition. If $a > 0$, it follows $\left(ab_b c_{n-b-1}\right)^2 = b_{b+1}c_{n-b-1}$. If $a = 0$ follows $\left(0b_{b-1}c_{n-b}\right)^2 = 0c_{n-1}$. So, in all the cases the interior of the triangle is not a semiring. \hfill $\Box$

\vspace{3mm}

When $n = 3$, the interior of the least triangle $\triangle^{(3)}\{0,1,2\}$  is one-element semiring and this element is an identity $\mathbf{i} = 123$.

\vspace{3mm}

Let $a, b, c, x, y, z \in \mathcal{C}_n$, $a < b < c$ and $x < y < z$. We consider the map
$$\Phi : \triangle^{(n)}\{a,b,c\} \rightarrow \triangle^{(n)}\{x,y,z\}$$
such that $\Phi\left(a_kb_\ell c_{n-k-\ell}\right) = x_ky_\ell z_{n-k-\ell}$, where $0 \leq k \leq n-1$, $0 \leq k \leq n-1$ and $0 \leq n - k -\ell \leq n-1$. Obviously, $\Phi$ is order-preserving. Hence, the additive semigroups of any two triangles $\triangle^{(n)}\{a,b,c\}$ and $\triangle^{(n)}\{x,y,z\}$ are isomorphic. But $\triangle^{(n)}\{a,b,c\}$ and $\triangle^{(n)}\{x,y,z\}$ are nonisomorphic  semirings.

\vspace{3mm}

\textbf{Proposition} \z \textsl{In arbitrary triangle $\triangle^{(n)}\{a,b,c\}$, where $n > 3$, there is at least one right identity and there are not any left identities.}

\emph{Proof.}
The least idempotent of $\mathcal{STR}^{(n)}\{a,b\}$ is  endomorphism $a_bb_{n-b}$ and the least idempotent of $\mathcal{STR}^{(n)}\{a,c\}$ is  endomorphism $a_cc_{n-c}$. Their sum is $\varepsilon = a_bb_{c-b}c_{n-c}$. Let $\alpha \in \mathcal{STR}^{(n)}\{a,b\}$. If $\alpha(a) = a$, or $\alpha(a) = b$, or $\alpha(a) = c$, then it follows $(\alpha\cdot \varepsilon)(a) = a$, or $(\alpha\cdot \varepsilon)(a) = b$, or $(\alpha\cdot \varepsilon)(a) = a$, respectively. The same is valid if we replace $a$ with $b$, or $a$ with $c$. So,  endomorphism $\varepsilon = a_bb_{c-b}c_{n-c}$ is a right identity of  triangle $\triangle^{(n)}\{a,b,c\}$.

If $b > a + 1$,  endomorphism $a_{a+1}b_{n-a-1}$ is another (in all the cases, the biggest) idempotent of $\mathcal{STR}^{(n)}\{a,b\}$. Now, in a similar way as above, it is easy to check that sum $a_{a+1}b_{n-a-1} + a_cc_{n-c} = a_{a+1}b_{c-a-1}c_{n-c}$ is another right identity of $\triangle^{(n)}\{a,b,c\}$.

Let us, similarly, suppose that $c > b + 1$. Then  endomorphism $a_{b+1}c_{n-b-1}$ is another, different from $a_cc_{n-c}$, idempotent of $\mathcal{STR}^{(n)}\{a,c\}$. Now  sum $a_bb_{n-b} + a_{b+1}c_{n-b-1} = a_bbc_{n-b-1}$ is a right identity of $\triangle^{(n)}\{a,b,c\}$.

If there are two right identities of the triangle, it implies that there is not a left identity in this semiring. So, we consider the case when there is only one idempotent of $\mathcal{STR}^{(n)}\{a,b\}$ and there is only one idempotent of $\mathcal{STR}^{(n)}\{a,c\}$. It is possible when $b = a + 1$ and $c = a + 2$. Thus, it is enough to prove that there is not a left identity in any triangle $\triangle^{(n)}\{a,a+1,a+2\}$. We consider two cases.

\vspace{2mm}

\emph{Case 1.} Let $a \geq 1$ and $\alpha = (a+1)(a+2)_{n-1}$. Then $\alpha(a) = \alpha(a+1) = \alpha(a+2) = a + 2$. Hence, for any
$\beta \in \triangle^{(n)}\{a,a+1,a+2\}$ we find $\beta\cdot \alpha = \overline{a + 2}$. Since $\beta\cdot \overline{a + 2} = \overline{a + 2}$ it follows that
$\alpha \, \sim_\ell \; \overline{a + 2}$. So, there are two left-similar elements of $\triangle^{(n)}\{a,a+1,a+2\}$ and hence there is not a left identity.

\vspace{2mm}

\emph{Case 2.} Let $a = 0$. In semiring $\triangle^{(n)}\{0,1,2\}$ we choose endomorphism $\alpha = 1_{n-1}2$. Then for any $\beta \in \triangle^{(n)}\{0,1,2\}$ it follows $\beta\cdot \alpha = \overline{1}$. Since $\beta\cdot \overline{1} = \overline{1}$, we show that $\alpha \, \sim_\ell \; \overline{1}$. So, there are two left-similar elements of $\triangle^{(n)}\{0,1,2\}$ and  there is not a left identity in this triangle. \hfill $\Box$

\vspace{3mm}

Since there is a right identity in semiring $\triangle^{(n)}\{a,b,c\}$, it follows

\vspace{3mm}

\textbf{Corollary} \z  \textsl{In an arbitrary triangle $\triangle^{(n)}\{a,b,c\}$, where $n > 3$, there are  not right-similar endomorphisms.}

\vspace{3mm}

Immediately it follows

\vspace{3mm}

\textbf{Corollary} \z \textsl{Only in the least triangle $\triangle^{(3)}\{0,1,2\}$ there is an identity $\mathbf{i} = 123$.}

\vspace{5mm}

\noindent{\large \bf 5. Layers in a Triangle}

\vspace{3mm}

Since any triangle $\triangle^{(n)}\{a,b,c\}$ is a 2--simplex, we define layers as in section 2. Any layer of a triangle is a chain. So, the elements of  layer
$\mathcal{L}^{k}_{a}\left(\triangle^{(n)}\{a,b,c\}\right)$, where $k = 0, \ldots, n-1$, are the following $n-k+1$ endomorphisms:
$$a_kb_{n-k} < a_kb_{n-k-1}c < \cdots < a_kc_{n-k}.$$

Similarly, we can represent the elements of $\mathcal{L}^{k}_{b}\left(\triangle^{(n)}\{a,b,c\}\right)$ and $\mathcal{L}^{k}_{c}\left(\triangle^{(n)}\{a,b,c\}\right)$, where $k = 0, \ldots, n-1$.

When the least and the biggest element of the layer are idempotents, we call this layer a {\emph{basic layer}}. Let us consider  layer $\mathcal{L}^{k}_{a}\left(\triangle^{(n)}\{a,b,c\}\right)$ with respect to  vertex $\overline{a}$.. Since idempotents of $\mathcal{STR}^{(n)}\{a,b\}$ are endomorphisms
$a_bb_{n-b}, \ldots, a_{a+1}b_{n-a-1}$ and, similarly, idempotents of $\mathcal{STR}^{(n)}\{a,c\}$ are $a_cc_{n-c}, \ldots, a_{a+1}c_{n-a-1}$,
it follows that $\mathcal{L}^{k}_{a}\left(\triangle^{(n)}\{a,b,c\}\right)$ is a basic layer only if ${k = a + 1, \ldots, b}$.
All the elements of the layer are endomorphisms $\alpha = a_kb_{n-k-i}c_i$, where ${i = 0, \ldots, n-k}$.
It is easy to see that if $i \leq n - c -1$, then $\alpha(b) = \alpha(c) = b$. Such endomorphisms are called {\emph{left elements}} of the layer. If $n -c \leq i \leq n - b -1$, then $\alpha(b) =b$ and $\alpha(c) = c$, so, $\alpha$ is an idempotent which is  a right identity of the triangle. If $i \geq n -b$, then $\alpha(b) = \alpha(c) = c$. Such endomorphisms are called {\emph{right elements}} of the layer. Let $\alpha$ and $\beta$ be left elements of the layer. If $\alpha(x) = b$, where $x \in \mathcal{C}_n$, then $(\alpha\cdot \beta)(x) = \beta(b) = b$. If $\alpha(x) = c$, where $x \in \mathcal{C}_n$, then $(\alpha\cdot \beta)(x) = \beta(c) = b$. Hence, $\alpha\cdot \beta = a_kb_{n-k}$. Similarly, if $\alpha$ and $\beta$ are right elements of the layer, then $\alpha\cdot \beta = a_kc_{n-k}$. Let $\alpha$ be a left element and $\beta$ be a right element of the layer. If $\alpha(x) = b$, where $x \in \mathcal{C}_n$, then $(\alpha\cdot \beta)(x) = \beta(b) = c$. If $\alpha(x) = c$, where $x \in \mathcal{C}_n$, then $(\alpha\cdot \beta)(x) = \beta(c) = c$. Hence, $\alpha\cdot \beta = a_kc_{n-k}$. Similarly, $\beta\cdot \alpha = a_kb_{n-k}$. Thus we obtain

\vspace{3mm}

\textbf{Proposition} \z \textsl{Any basic layer $\mathcal{L}^{k}_{a}\left(\triangle^{(n)}\{a,b,c\}\right)$, $k = a + 1, \ldots, b$, is a semi\-ring.}

\vspace{3mm}

From the proof of last proposition it follows that in  layer $\mathcal{L}^{k}_{a}\left(\triangle^{(n)}\{a,b,c\}\right)$ the first $n-c$ elements $\alpha$ are left elements, i.e. $\alpha^2 = a_kb_{n-k}$, the next $c-b$ endomorphisms are idempotents and the last $b - k +1$ elements $\alpha$ are right elements, i.e. $\alpha^2 = a_kc_{n-k}$. If there is a string $\mathcal{STR}^{(n_0)}\{a_0,b_0\}$ with $n-c$ $\,a_0$--nilpotent elements, $c -b$ idempotents and $b - k +1$ $\,b_0$--nilpotent elements are $b-k+1$, then two semirings $\mathcal{L}^{k}_{a}\left(\triangle^{(n)}\{a,b,c\}\right)$ and $\mathcal{STR}^{(n_0)}\{a_0,b_0\}$ will be isomorphic. From  system $\left|\begin{array}{l} n_0 - b_0 = n -c\\ b_0 - a_0 = c - b\\a_0 + 1 = b - k + 1 \end{array}\right.$ we find $n_0 = n-k$, $a_0 = b -k$ and $b_0 = c - k$. So, we prove

\vspace{3mm}

\textbf{Proposition} \z \textsl{For any $n \geq 3$, $a, b \in \mathcal{C}_n$, $a < b$ and $k = a+1, \ldots, b$, semirings $\mathcal{L}^{k}_{a}\left(\triangle^{(n)}\{a,b,c\}\right)$ and $\mathcal{STR}^{(n-k)}\{b-k,c-k\}$ are isomorphic.}

\vspace{3mm}

The basic layers $\mathcal{L}^{k}_{b}\left(\triangle^{(n)}\{a,b,c\}\right)$ are not closed under the multiplication in the general case. For instance, see fig,1, where the layer $\mathcal{L}^{2}_{2}\left(\triangle^{(4)}\{1,2,3\}\right)$ is a basic layer, but for his middle element $12_23$ we find $(12_23)^2 = 2_33 \notin \mathcal{L}^{2}_{2}\left(\triangle^{(4)}\{1,2,3\}\right)$.

Now let us consider  layer $\mathcal{L}^{k}_{c}\left(\triangle^{(n)}\{a,b,c\}\right)$ with respect to  vertex $\overline{c}$. Since idempotents of $\mathcal{STR}^{(n)}\{a,c\}$ are $a_cc_{n-c}, \ldots, a_{a+1}c_{n-a-1}$ and idempotents of $\mathcal{STR}^{(n)}\{b,c\}$ are $b_cc_{n-c}, \ldots, b_{b+1}c_{n-b-1}$, it follows that $\mathcal{L}^{k}_{c}\left(\triangle^{(n)}\{a,b,c\}\right)$ is a basic layer only if ${k = n - c, \ldots, n - b -1}$.
All the elements of the layer are endomorphisms $\alpha = a_ib_{n-k-i}c_k$, where ${i = 0, \ldots, n-k}$. If $b+1 \leq i \leq n-k$, then $\alpha(a) = \alpha(b) = a$. We call these endomorphisms (as in the previous case) left elements of the layer. When $a+ 1 \leq i \leq b$, it follows $\alpha(a) = a$, $\alpha(b) = b$, that is $\alpha$ is an idempotent which is a right identity of the triangle. If $0 \leq i \leq a$, then $\alpha(a) = \alpha(b) = b$. These endomorphisms are the right elements of the layer.  By the same way, as for the basic layers with respect to vertex $\overline{a}$, we prove here that:

1. If $\alpha$ and $\beta$ are left elements of the layer, then $\alpha\cdot \beta = a_{n-k}c_{k}$.

2. If $\alpha$ and $\beta$ are right elements of the layer, then $\alpha\cdot \beta = b_{n-k}c_{k}$.

3. If $\alpha$ is a left but $\beta$ is a right element of the layer, then $\alpha\cdot \beta = b_{n-k}c_{k}$ and $\beta\cdot \alpha = a_{n-k}c_{k}$

So, we obtain

\vspace{3mm}

\textbf{Proposition} \z \textsl{Any basic layer $\mathcal{L}^{k}_{c}\left(\triangle^{(n)}\{a,b,c\}\right)$, where $k = n - c, \ldots, n-b-1$, is a semi\-ring.}

\vspace{3mm}

Now we search a string $\mathcal{STR}^{(n_0)}\{a_0,b_0\}$ with $n-k-b$ $\;a_o$--nilpotent elements, $b-a$ idempotents and $a + 1$ $\;b_0$--nilpotent elements. It is easy to find that $n_0 = n -k$, $a_0 = a$ and $b_0 = b$. Thus we prove

\vspace{3mm}

\textbf{Proposition} \z \textsl{For any $n \geq 3$, $b, c \in \mathcal{C}_n$, $b < c$ and $k = n-c, \ldots, n-b-1$,  semirings $\mathcal{L}^{k}_{c}\left(\triangle^{(n)}\{a,b,c\}\right)$ and $\mathcal{STR}^{(n-k)}\{a,b\}$ are isomorphic.}

\vspace{3mm}

In  triangle $\triangle^{(4)}\{1,2,3\}$ (fig. 1) we obtain that there are many left-similar endomorphisms. In order to prove that there are such elements in any triangle $\triangle^{(n)}\{a,b,c\}$ we first consider the case $a > 0$. Then $ac_{n-1} \sim_\ell \overline{c}$. Actually, if $\alpha = ac_{n-1}$, easily follows that $\alpha(a) = \alpha(b) = \alpha(c) = c$. Thus, for any $\beta \in \triangle^{(n)}\{a,b,c\}$ we have $\beta\cdot \alpha = \overline{c} = \beta\cdot \overline{c}$, that is $\alpha \sim_\ell \overline{c}$. Note that $\alpha$ and $\overline{c}$ are not right identities of the triangle.

Let $a = 0$. In $\triangle^{(n)}\{0,b,c\}$ we consider the biggest basic layer $\mathcal{L}^{1}_{0}\left(\triangle^{(n)}\{0,b,c\}\right)$. The elements of this layer are $0b_{n-1} < \cdots < 0bc_{n-2} < 0c_{n-1}$. Now, if $b > 1$, it follows that $0bc_{n-2} \sim_\ell 0c_{n-1}$. Indeed, let $\alpha = 0bc_{n-1}$. Since $\alpha(0) = 0$, $\alpha(b) = \alpha(c) = c$, it follows that if $\beta \in \triangle^{(n)}\{0,b,c\}$, then $\beta\cdot \alpha = \beta\cdot 0c_{n-1}$, i.e. $\alpha \sim_\ell 0c_{n-1}$. Note that $\alpha$ and $0c_{n-1}$ are not right identities of the triangle.

Now let us consider  triangle $\triangle^{(n)}\{0,1,c\}$ and endomorphisms $01_{n-1}$ and $01_{n-2}$. Let $c < n-1$. Then for $\alpha = 01_{n-2}c$ we have $\alpha(0) = 0$, $\alpha(b) = \alpha(c) = 1$. So, for any $\beta \in \triangle^{(n)}\{0,1,c\}$, it follows $\beta\cdot \alpha = \beta\cdot 01_{n-1}$, that is $01_{n-1}c \sim_\ell 01_{n-1}$. Note that these endomorphisms are not right identities.

 Finally let $c = n-1$. We consider  $\triangle^{(n)}\{0,1,n-1\}$, where $n > 3$. Note, that in the least triangle $\triangle^{(3)}\{0,1,2\}$ there are not left-similar elements since there is an identity. Now let us consider  endomorphisms $\alpha = 0_{n-1}1$ and $\beta = 0_{n-2}1_2$. We have
 $\alpha(0) = \beta(0) = 0$,  $\alpha(1) = \beta(1) = 0$,  $\alpha(n-1) = \beta(n-1) = 1$. So, for any $\gamma \in \triangle^{(n)}\{0,1,n-1\}$, it follows $\gamma\cdot \alpha = \gamma\cdot \beta$, i.e. $\alpha \sim_\ell \beta$. Note that $\alpha$ and $\beta$ are not right identities. Hence, we prove

\vspace{3mm}

\textbf{Proposition} \z \textsl{For any $n > 3$ the triangle $\triangle^{(n)}\{a,b,c\}$ contains a pair of elements which are left-similar endomorphisms and are not right identities.}

\vspace{5mm}

\noindent{\large \bf 6. Idempotents and Nilpotent Elements of a Triangle}

\vspace{3mm}

By a {\emph{boundary idempotent}} of  triangle $\triangle^{(n)}\{a,b,c\}$ we understand an idempotent of any of the strings of this triangle. The idempotent of the interior of the triangle is called {\emph{interior idempotent}}. From Proposition 17, it follows that the set of all boundary idempotents is closed under the multiplication.

\vspace{3mm}

\textbf{Proposition} \z \textsl{The interior idempotents of the triangle $\triangle^{(n)}\{a,b,c\}$ are just the right identities of the triangle.}

\emph{Proof.} As we know, [2],  endomorphism $\alpha \in \widehat{\mathcal{E}}_{\mathcal{C}_n}$ with $s$ fixed points $k_1, \ldots, k_s$, ${1 \leq s \leq
n-1}$, is an idempotent if and only if  $\im(\alpha) = \{k_1, \ldots, k_s\}$. So, for any interior idempotent of $\triangle^{(n)}\{a,b,c\}$, it follows that $a$, $b$ and $c$ are fixed points of $\alpha$. Then, for any $\beta \in \triangle^{(n)}\{a,b,c\}$ and $x \in \mathcal{C}_n$ easily follows $(\beta\cdot \alpha)(x) = \alpha(\beta(x)) = \beta(x)$. Hence, $\alpha$ is a right identity of $\triangle^{(n)}\{a,b,c\}$.

Conversely, if $\alpha$ is a right identity, then, obviously, $\alpha$ is an idempotent. If we assume that $\alpha$ is a boundary idempotent, say $\alpha = a_kb_{n-k}$, where $k = a+1, \ldots, b$, then $\overline{c}\cdot a_kb_{n-k} = \overline{b}$. This contradicts the choice that $\alpha$ is a right identity. \hfill $\Box$

\vspace{3mm}

The set of right identities of $\triangle^{(n)}\{a,b,c\}$ is denoted by $\mathcal{RI}\left(\triangle^{(n)}\{a,b,c\}\right)$.

\vspace{3mm}

 For any triangle $\triangle^{(n)}\{a,b,c\}$  endomorphism $a_{a+1}c_{n-a-1}$ is an idempotent of $\mathcal{STR}^{(n)}\{a,c\}$ and  endomorphism $b_{b+1}c_{n-b-1}$ is an idempotent of $\mathcal{STR}^{(n)}\{b,c\}$. Now, it follows that $a_{a+1}c_{n-a-1} + b_{b+1}c_{n-b-1} = b_{a+1}c_{n-a-1} \in N^{[a]}\left(\mathcal{STR}^{(n)}\{a,b\}\right)$ since $(b_{a+1}c_{n-a-1})^2 = \overline{c}$. So, the set of idempotents of any triangle is not a semiring. But this does not mean that adding  idempotents is a ``bad idea''.

Any discrete neighborhood of a vertex of the triangle $\triangle^{(n)}\{a,b,c\}$ (see fig. 2) can be represented as a triangle in a geometrical sense. So, such subset of $\triangle^{(n)}\{a,b,c\}$ is called {\emph{geometric triangle}}. In fig. 2 the endomorphisms $a_{n-m}c_m$, $b_{n-m}c_m$ and $\overline{c}$ are the ``vertices'' and subsets of $\mathcal{STR}^{(n)}\{a,c\}$ and $\mathcal{STR}^{(n)}\{b,c\}$ consisting of $n-m+1$ endomorphisms and also the layer $\mathcal{L}^{n-m}_{c}\left(\triangle^{(n)}\{a,b,c\}\right)$ are the ``sides'' of this geometric triangle.
 Similarly, geometric triangles are called subsets depicted in figures 3 and 4.

\begin{figure}[h]\centering
  \includegraphics[width=150mm]{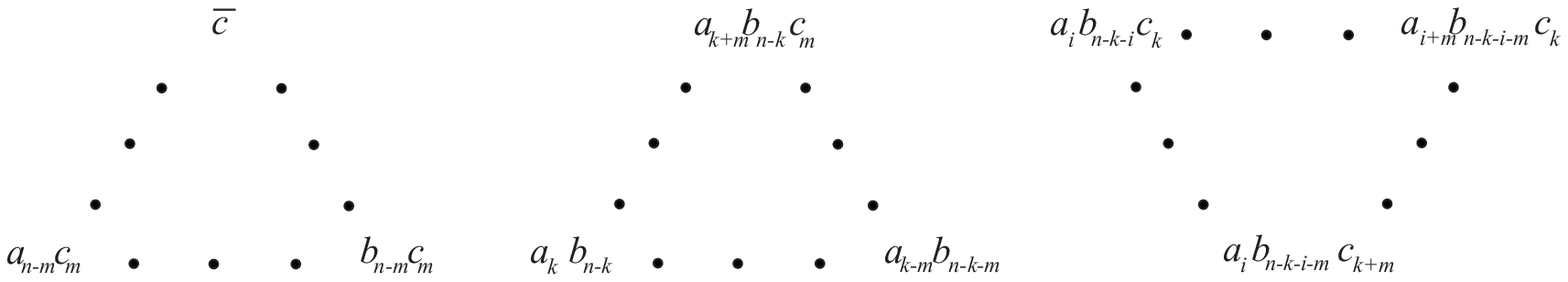}\\
\end{figure}

\centerline{\small \hphantom{aaaa} Figure 2 \hphantom{aaaaaaaaaaaaaaaaaa} Figure 3 \hphantom{aaaaaaaaaaaaaaaaaaaa} Figure 4 \hphantom{aaaaa}}

\vspace{3mm}

Geometric triangles are, in general, not semirings. But some of them are semirings, for example,   geometric triangles
$\mathcal{DN}^{n-1}_a = \{\overline{a}, a_{n-1}b, a_{n-1}c\}$, $\mathcal{DN}^{n-1}_b = \{\overline{b}, ab_{n-1}, b_{n-1}c\}$ and $\mathcal{DN}^{n-1}_c = \{\overline{c}, ac_{n-1}, bc_{n-1}\}$ are subsemirings of $\triangle^{(n)}\{a,b,c\}$ from Proposition 2.

From Theorem 5 it follows that the discrete neighborhood of  vertex $\overline{a}$ of $\triangle^{(n)}\{a,b,c\}$ containing the biggest layer $\mathcal{L}^{a+1}_{a}\left(\triangle^{(n)}\{a,b,c\}\right)$ is  semiring $\mathcal{DN}^{n-a-1}_a = \triangle^{(n)}\{a,b,c\}\cap {\mathcal{E}}^{(a)}_{\mathcal{C}_n}$. This semiring is a geometric triangle whose ``vertices'' the idempotent endomorphisms $\overline{a}$, $a_{a+1}b_{n-a-1}$ and $a_{a+1}c_{n-a-1}$. The ``sides'' of this geometric triangle are also semirings, since $\mathcal{L}^{a+1}_{a}\left(\triangle^{(n)}\{a,b,c\}\right)$ is a semiring (see Proposition 27) and other ``sides'' are semirings of type $A_{a+1}$ (see the end of section 3).

Similarly, from Theorem 5 we know that the discrete neighborhood of the vertex $\overline{c}$ of $\triangle^{(n)}\{a,b,c\}$ containing the biggest layer $\mathcal{L}^{n-c}_{c}\left(\triangle^{(n)}\{a,b,c\}\right)$ is  semiring
 $\mathcal{DN}^{c}_c = \triangle^{(n)}\{a,b,c\}\cap {\mathcal{E}}^{(c)}_{\mathcal{C}_n}$.  This semiring is also a geometric triangle whose ``vertices'' are the idempotent endomorphisms $a_cc_{n-c}$, $b_cc_{n-c}$ and $\overline{c}$.

The intersection $\mathcal{DN}^{n-a-1}_a\cap \mathcal{DN}^{c}_c$ is a semiring consisting of all the endomorphisms with fixed points $a$ and $b$. Thus we construct  a new geometric triangle whose ``vertices'' are the idempotent endomorphisms $a_cc_{n-c}$,  $a_{a+1}b_{c-a-1}c_{n-c}$ and $a_{a+1}c_{n-a-1}$. This semiring is called an {\emph{idempotent triangle}} of $\triangle^{(n)}\{a,b,c\}$ and is denoted by $\mathcal{IT}\left(\triangle^{(n)}\{a,b,c\}\right)$.

\vspace{3mm}

\textbf{Theorem} \z \textsl{For any triangle $\triangle^{(n)}\{a,b,c\}$, $n \geq 3$, the set of right identities $\mathcal{RI}\left(\triangle^{(n)}\{a,b,c\}\right)$ is a subsemiring of $\mathcal{IT}\left(\triangle^{(n)}\{a,b,c\}\right)$ of order $(b-a)(c-a)$. The set $\mathcal{IT}\left(\triangle^{(n)}\{a,b,c\}\right)\backslash \mathcal{RI}\left(\triangle^{(n)}\{a,b,c\}\right)$ is a subsemiring of $\mathcal{IT}\left(\triangle^{(n)}\{a,b,c\}\right)$ of order $ \frac{1}{2}((c-b)^2 + (b-a)^2 + c- a)$. The semirings $Id\left(\mathcal{STR}^{(n)}\{a,c\}\right)$ and $\mathcal{IT}\left(\triangle^{(n)}\{a,b,c\}\right)\backslash \mathcal{RI}\left(\triangle^{(n)}\{a,b,c\}\right)$ are ideals of $\mathcal{IT}\left(\triangle^{(n)}\{a,b,c\}\right)$.}

\emph{Proof}. Let $\alpha = a_kb_{n-k-j}c_j$ be a right identity. From the last proposition it follows that $a$, $b$ and $c$ are the fixed points of $\alpha$. Since $\alpha(a) = a$, it follows $k \geq a + 1$. Since $\alpha(c) = c$, we have $j \geq n-c$. Finally, since $\alpha(b) = b$, it follows $a+1 \leq k \leq b$, $n-c \leq j \leq n - a - 1$ and $n - k - j \geq 1$. Hence, $a_kb_{n-k-j}c_j = a_kb_{n-k} + a_{n-j}c_j$, where $a_kb_{n-k} \in Id\left(\mathcal{STR}^{(n)}\{a,b\}\right)$ and $a_{n-j}c_{j} \in Id\left(\mathcal{STR}^{(n)}\{a,c\}\right)$. Actually we showed that $$a_kb_{n-k-j}c_j = \mathcal{L}^{k}_{a}\left(\triangle^{(n)}\{a,b,c\}\right)\cap \mathcal{L}^{\ell}_{c}\left(\triangle^{(n)}\{a,b,c\}\right),$$
 where $k = a+1, \ldots, b$ and
$\ell = n - c, \ldots, n-b-1$. Thus we prove that any right identity is an intersection of basic layer with respect to $\overline{a}$ and basic layer with respect to $\overline{c}$. Since there are $b-a$ basic layers with respect to $\overline{a}$ and $c-b$ basic layers with respect to $\overline{c}$, it follows that all right identities are $(b-a)(c-b)$. Since all elements of $\mathcal{RI}\left(\triangle^{(n)}\{a,b,c\}\right)$ are right identities of $\triangle^{(n)}\{a,b,c\}$, it follows that the set $\mathcal{RI}\left(\triangle^{(n)}\{a,b,c\}\right)$ is closed under the multiplication.

Let $\alpha \in \mathcal{L}^{k_1}_{a}\left(\triangle^{(n)}\{a,b,c\}\right)$, $\beta \in \mathcal{L}^{k_2}_{a}\left(\triangle^{(n)}\{a,b,c\}\right)$ and $\alpha, \beta \in \mathcal{L}^{\ell}_{c}\left(\triangle^{(n)}\{a,b,c\}\right)$, where  where $k = a+1, \ldots, b$ and $\ell = n - c, \ldots, n-b-1$. If we assume that $k_1 \leq k_2$, then $\alpha + \beta = \alpha$. Similarly, if $\alpha$ and $\beta$ are endomorphisms of the same layer with respect to $\overline{a}$ and $\alpha \in \mathcal{L}^{\ell_1}_{c}\left(\triangle^{(n)}\{a,b,c\}\right)$, $\beta \in \mathcal{L}^{\ell_2}_{c}\left(\triangle^{(n)}\{a,b,c\}\right)$, where $\ell_1 \geq \ell_2$, then $\alpha + \beta = \alpha$. Finally, let
 $\alpha = \mathcal{L}^{k_1}_{a}\left(\triangle^{(n)}\{a,b,c\}\right)\cap \mathcal{L}^{\ell_1}_{c}\left(\triangle^{(n)}\{a,b,c\}\right)$ and
$\beta = \mathcal{L}^{k_2}_{a}\left(\triangle^{(n)}\{a,b,c\}\right)\cap \mathcal{L}^{\ell_2}_{c}\left(\triangle^{(n)}\{a,b,c\}\right)$, where $k_1 \leq k_2$ and $\ell_1 \leq \ell_2$. Then we take $\gamma = \mathcal{L}^{k_1}_{a}\left(\triangle^{(n)}\{a,b,c\}\right)\cap \mathcal{L}^{\ell_2}_{c}\left(\triangle^{(n)}\{a,b,c\}\right)$ and $\alpha + \beta = \gamma$. Hence, $\mathcal{RI}\left(\triangle^{(n)}\{a,b,c\}\right)$ is a subsemiring of $\mathcal{IT}\left(\triangle^{(n)}\{a,b,c\}\right)$ of order $(b-a)(c-a)$.

Since $\mathcal{IT}\left(\triangle^{(n)}\{a,b,c\}\right) \subset \mathcal{DN}^{n-a-1}_a$, it follows that there are not any $b$--nilpotent and $c$--nilpotent endomorphisms of $\mathcal{IT}\left(\triangle^{(n)}\{a,b,c\}\right)$.
Since $\mathcal{IT}\left(\triangle^{(n)}\{a,b,c\}\right) \subset \mathcal{DN}^{c}_c$ it follows that there are not any $a$--nilpotent endomorphisms of $\mathcal{IT}\left(\triangle^{(n)}\{a,b,c\}\right)$. But we prove in Proposition 29 that all the left elements from some basic layer with respect to $\overline{c}$ are square roots of endomorphisms of $Id\left(\mathcal{STR}^{(n)}\{a,c\}\right)$. So, the idempotent triangle consist of all the right identities, all the elements of $Id\left(\mathcal{STR}^{(n)}\{a,c\}\right)$ and all the square roots of elements of $Id\left(\mathcal{STR}^{(n)}\{a,c\}\right)$. The vertices of the idempotent triangle are $a_cc_{n-c}$,  $a_{a+1}b_{c-a-1}c_{n-c}$ and $a_{a+1}c_{n-a-1}$. So, any ``side'' of this geometric triangle consists of $c-a$ elements. Then there are $\frac{1}{2}(c-a)(c-a+1)$ endomorphisms of $\mathcal{IT}\left(\triangle^{(n)}\{a,b,c\}\right)$. Thus it follows that the elements of set $\mathcal{IT}\left(\triangle^{(n)}\{a,b,c\}\right)\backslash \mathcal{RI}\left(\triangle^{(n)}\{a,b,c\}\right)$ are\\ $\frac{1}{2}(c-a)(c-a+1) - (c-b)(b-a) = \frac{1}{2}((c-b)^2 + (b-a)^2 + c- a)$.

Now we consider a partition of  idempotent triangle $\mathcal{IT}\left(\triangle^{(n)}\{a,b,c\}\right)$ into three parts. The first one is a geometric triangle with ``vertices'' $a_cc_{n-c}$, $a_{b+1}b_{c-b-1}c_{n-c}$ and $a_{b+1}c_{n-b-1}$. This triangle is a subset of $\mathcal{DN}^{c}_c$ and there are not any common elements of triangle and the basic layers with respect to $\overline{a}$. Since all the endomorphisms of the triangle with an  exception of elements $Id\left(\mathcal{STR}^{(n)}\{a,c\}\right)$ are left elements of the basic layers with respect to $\overline{c}$, the triangle is denoted by $L_{\triangle}$.
The second part of  $\mathcal{IT}\left(\triangle^{(n)}\{a,b,c\}\right)$ is a geometric triangle with ``vertices'' $a_bc_{n-b}$, $a_{a+1}b_{b-a-1}c_{n-b}$ and $a_{a+1}c_{n-a-1}$. This triangle is a subset of $\mathcal{DN}^{n-a-1}_a$ and there are not any common elements of the triangle and the basic layers with respect to $\overline{c}$. Since all the endomorphisms of the triangle with an  exception of elements $Id\left(\mathcal{STR}^{(n)}\{a,c\}\right)$ are right elements of the basic layers with respect to $\overline{a}$, the triangle is denoted by $R_{\triangle}$. The third part of the idempotent triangle is  semiring $\mathcal{RI}\left(\triangle^{(n)}\{a,b,c\}\right)$ whose elements are the intersections of all the basic layers with respect to $\overline{a}$ and  $\overline{c}$. So, $L_{\triangle}\cup R_{\triangle} = \mathcal{IT}\left(\triangle^{(n)}\{a,b,c\}\right)\backslash \mathcal{RI}\left(\triangle^{(n)}\{a,b,c\}\right)$.

Let $\alpha \in L_{\triangle}$ and $\beta \in R_{\triangle}$. Then $\alpha = a_mb_{\ell-m}c_{n-\ell}$, where $\ell = b + 1, \ldots, c$, $\ell - m \geq 0$ and $\alpha(a) = \alpha(b) = a$. Similarly, $\beta = a_mb_{\ell-m}c_{n-\ell}$, where $m = a+1, \ldots, b$,  $\ell - m \geq 0$ and $\beta(b) = \beta(c) = c$.

Let $a_kc_{n-k} \in Id\left(\mathcal{STR}^{(n)}\{a,c\}\right)$. Then $k = a+1, \ldots, c$. We find:

1. $a_kc_{n-k}\cdot \alpha = a_kc_{n-k}$ for $\alpha \in L_{\triangle}$.

2. $a_kc_{n-k}\cdot \beta = a_kc_{n-k}$ for $\beta \in R_{\triangle}$.

3. $a_kc_{n-k}\cdot \varepsilon = a_kc_{n-k}$ for $\varepsilon \in \mathcal{RI}\left(\triangle^{(n)}\{a,b,c\}\right)$.

So, the elements of  semiring $Id\left(\mathcal{STR}^{(n)}\{a,c\}\right)$ are left zeroes of semiring $\mathcal{IT}\left(\triangle^{(n)}\{a,b,c\}\right)$.

Now we calculate:

4. $\alpha\cdot a_kc_{n-k} = a_mb_{\ell-m}c_{n-\ell}\cdot a_kc_{n-k} = a_{\ell}c_{n-\ell} \in Id\left(\mathcal{STR}^{(n)}\{a,c\}\right)$.

5. $\beta\cdot a_kc_{n-k} = a_mb_{\ell-m}c_{n-\ell}\cdot a_kc_{n-k} = a_{m}c_{n-m} \in Id\left(\mathcal{STR}^{(n)}\{a,c\}\right)$.

6. For any $\varepsilon = a_tb_{n-t-j}c_j \in \mathcal{RI}\left(\triangle^{(n)}\{a,b,c\}\right)$, where $t = a+1, \ldots, b$ and\\ $j = n-c, \ldots, n-a-1$ it follows
$\varepsilon\cdot a_kc_{n-k} = \left\{ \begin{array}{ll} a_tc_{n-t}& \mbox{for}\;\; k = a+1, \ldots, b\\ a_{n-j}c_j & \mbox{for}\; \; k = b+1, \ldots, c \end{array} \right. $. Thus, in all cases $\varepsilon\cdot a_kc_{n-k} \in Id\left(\mathcal{STR}^{(n)}\{a,c\}\right)$. Hence, the semiring $Id\left(\mathcal{STR}^{(n)}\{a,c\}\right)$ is an ideal of $\mathcal{IT}\left(\triangle^{(n)}\{a,b,c\}\right)$.

In order to prove that $\mathcal{IT}\left(\triangle^{(n)}\{a,b,c\}\right)\backslash \mathcal{RI}\left(\triangle^{(n)}\{a,b,c\}\right)$ is an  ideal of $\mathcal{IT}\left(\triangle^{(n)}\{a,b,c\}\right)$ we find:

7. For any $\alpha = a_mb_{\ell-m}c_{n-\ell} \in  L_{\triangle}$ and $\alpha_1 = a_{m_1}b_{{\ell_1}-m_1}c_{n-{\ell_1}} \in  L_{\triangle}$, where $\ell, \ell_1 = b + 1, \ldots, c$, it follows $\alpha\cdot \alpha_1 = a_{\ell}c_{n-\ell} \in Id\left(\mathcal{STR}^{(n)}\{a,c\}\right)$.

8. For any $\beta = a_mb_{\ell-m}c_{n-\ell} \in  R_{\triangle}$ and $\beta_1 = a_{m_1}b_{{\ell_1}-m_1}c_{n-{\ell_1}} \in  R_{\triangle}$, where
 $m, m_1 = a + 1, \ldots, b$, it follows $\beta\cdot \beta_1 = a_{m}c_{n-m} \in Id\left(\mathcal{STR}^{(n)}\{a,c\}\right)$.

9. For any $\alpha = a_mb_{\ell-m}c_{n-\ell} \in  L_{\triangle}$ and $\beta = a_{m_1}b_{{\ell_1}-m_1}c_{n-{\ell_1}} \in  r_{\triangle}$, where
 $\ell = b + 1, \ldots, c$ and $m = a+1, \ldots, b$, it follows $\alpha\cdot \beta = a_{m}c_{n-m} \in Id\left(\mathcal{STR}^{(n)}\{a,c\}\right)$.

10. For any $\alpha = a_mb_{\ell-m}c_{n-\ell} \in  L_{\triangle}$ and $\beta = a_{m_1}b_{{\ell_1}-m_1}c_{n-{\ell_1}} \in  r_{\triangle}$, where $\ell = b + 1, \ldots, c$ and $m = a+1, \ldots, b$, it follows $\beta\cdot \alpha = a_{\ell}c_{n-\ell} \in Id\left(\mathcal{STR}^{(n)}\{a,c\}\right)$.

11. For any $\alpha  \in  L_{\triangle}$, $\beta \in  R_{\triangle}$ and $\varepsilon  \in \mathcal{RI}\left(\triangle^{(n)}\{a,b,c\}\right)$ it follows $\alpha\cdot \varepsilon = \alpha$ and $\beta\cdot \varepsilon = \beta$.

12. For any $\varepsilon = a_tb_{n-t-j}c_j \in \mathcal{RI}\left(\triangle^{(n)}\{a,b,c\}\right)$, where $t = a+1, \ldots, b$ and\\ ${j = n-c, \ldots, n-a-1}$, $\alpha = a_kb_{\ell-k}c_{n-\ell} \in  L_{\triangle}$, where $\ell = b + 1, \ldots, c$ and\\ $\beta = a_mb_{\ell-m}c_{n-\ell} \in  R_{\triangle}$, where $m = a + 1, \ldots, b$, it follows
$$\varepsilon\cdot \alpha = a_tb_{n-t-j}c_j\cdot a_kb_{\ell-k}c_{n-\ell} = a_{n-j}c_j \in Id\left(\mathcal{STR}^{(n)}\{a,c\}\right)\; \mbox{and}\;$$
$$\varepsilon\cdot \beta = a_tb_{n-t-j}c_j\cdot a_mb_{\ell-m}c_{n-\ell} = a_tc_{n-t} \in Id\left(\mathcal{STR}^{(n)}\{a,c\}\right).$$

The endomorphisms $\alpha$ of $L_{\triangle}$ are characterized in  triangle $\triangle^{(n)}\{a,b,c\}$ by equalities: $\alpha(a) = \alpha(b) = a$, $\alpha(c) = c$. So, if $\alpha, \beta \in L_{\triangle}$, then $(\alpha + \beta)(a) = (\alpha + \beta)(b) = a$ and $(\alpha + \beta)(c) = c$, i.e. $\alpha + \beta \in L_{\triangle}$. Similar reasonings we can have for triangle $R_{\triangle}$.

The biggest endomorphism of $L_{\triangle}$ is the vertex of this triangle $a_{b+1}c_{n-b-1}$. The least endomorphism of $R_{\triangle}$ is the vertex of the triangle $a_bc_{n-b}$ So, for any $\alpha \in L_{\triangle}$ and $\beta \in R_{\triangle}$ it follows $\alpha < a_{b+1}c_{n-b-1} < a_bc_{n-b} < \beta$. Hence, $\mathcal{IT}\left(\triangle^{(n)}\{a,b,c\}\right)\backslash \mathcal{RI}\left(\triangle^{(n)}\{a,b,c\}\right)$ is an  ideal of $\mathcal{IT}\left(\triangle^{(n)}\{a,b,c\}\right)$. \hfill $\Box$

\vspace{3mm}

From the proof of the last theorem it follows

\vspace{3mm}

\textbf{Corollary} \z  \textsl{The geometric triangles $L_{\triangle}$ and $R_{\triangle}$ are subsemirings of  triangle $\triangle^{(n)}\{a,b,c\}$.}

\vspace{3mm}

An immediate consequence of fact above is

\vspace{3mm}

\textbf{Corollary} \z  \textsl{For any triangle $\triangle^{(n)}\{a,b,c\}$, $n \geq 3$, the idempotent triangle is a disjoint union of subsemirings
$L_{\triangle}$, $R_{\triangle}$ and $\mathcal{RI}\left(\triangle^{(n)}\{a,b,c\}\right)$.}

\vspace{3mm}

Similarly to geometric triangles, we can consider  {\emph{geometric parallelograms}} and {\emph{geometric trapezoids}}. For example,  semiring $\mathcal{RI}\left(\triangle^{(n)}\{a,b,c\}\right)$ can be represented as a geometric parallelogram whose ``vertices'' are endomorphisms $a_bb_{c-b}c_{n-c}$, $a_{a+1}b_{c-a-1}c_{n-c}$, $a_{a+1}b_{b-a}c_{n-b-1}$ and $a_bbc_{n-b-1}$. (Note that exactly the last endomorphism is a boundary between the triangles $L_{\triangle}$ and $R_{\triangle}$.) The ``sides'' of this parallelogram are the idempotent parts of basic layers $\mathcal{L}^{n-c}_{c}\left(\triangle^{(n)}\{a,b,c\}\right)$, $\mathcal{L}^{n-b-1}_{c}\left(\triangle^{(n)}\{a,b,c\}\right)$, $\mathcal{L}^{a+1}_{a}\left(\triangle^{(n)}\{a,b,c\}\right)$ and $\mathcal{L}^{b}_{a}\left(\triangle^{(n)}\{a,b,c\}\right)$.

\vspace{3mm}

Now we consider the $a$--nilpotent elements of  triangle $\triangle^{(n)}\{a,b,c\}$. The set of all $a$--nilpotent elements of this triangle is denoted by $N^{[a]}\left( \triangle^{(n)}\{a,b,c\}\right)$. Since $N^{[a]}\left( \triangle^{(n)}\{a,b,c\}\right) = N^{[a]}_n\cap \triangle^{(n)}\{a,b,c\}$, similarly to Proposition 9, it follows

\vspace{3mm}

\textbf{Proposition} \z  \textsl{The set $N^{[a]}\left( \triangle^{(n)}\{a,b,c\}\right)$ is a subsemiring of $\triangle^{(n)}\{a,b,c\}$.}

\vspace{3mm}

The semiring $N^{[a]}\left( \triangle^{(n)}\{a,b,c\}\right)$ can be represented as a geometric trapezoid whose ``vertices'' are endomorphisms $\overline{a}$, $a_{b+1}b_{n-b-1}$, $a_{b+1}b_{c-b}c_{n-c-1}$ and $a_{c+1}c_{n-c-1}$ and whose ``sides'' are semiring $N^{[a]}\left(\mathcal{STR}^{(n)}\{a,b\}\right)$, the subsetset of  $n-c$ endomorphisms $\alpha$ from the left part of $\mathcal{L}^{b+1}_{a}\left(\triangle^{(n)}\{a,b,c\}\right)$ such that $\alpha(b) = a, \alpha(c) = b$,  the subsetset of  $c-b+1$ endomorphisms $\beta$ from the left part of $\mathcal{L}^{n-c-1}_{c}\left(\triangle^{(n)}\{a,b,c\}\right)$ such that $\alpha(b) = a, \alpha(c) = b$ and  semiring $N^{[a]}\left(\mathcal{STR}^{(n)}\{a,c\}\right)$. We find, as in the proof of the last theorem, that the order of this semiring is  $\left|N^{[a]}\left( \triangle^{(n)}\{a,b,c\}\right)\right| = \frac{1}{2}(n-c)(n + c - 2b + 1)$.

In the same way we can construct the following semirings:\\ $N^{[b]}\left( \triangle^{(n)}\{a,b,c\}\right) = N^{[b]}_n\cap \triangle^{(n)}\{a,b,c\}$ and $N^{[c]}\left( \triangle^{(n)}\{a,b,c\}\right) = N^{[c]}_n\cap \triangle^{(n)}\{a,b,c\}$.

The semiring $N^{[b]}\left( \triangle^{(n)}\{a,b,c\}\right)$ can be represented as a geometric parallelogram whose ``vertices'' are endomorphisms $a_ab_{n-a}$, $\overline{b}$, $b_{c+1}c_{n-c-1}$ and $a_ab_{c-a+1}c_{n-c-1}$ and whose ``sides'' are semiring  $N^{[b]}\left(\mathcal{STR}^{(n)}\{a,b\}\right)$,  semiring $N^{[b]}\left(\mathcal{STR}^{(n)}\{b,c\}\right)$,  the subsetset of  $a+1$ endomorphisms $\alpha$ from the right part of $\mathcal{L}^{n-c-1}_{c}\left(\triangle^{(n)}\{a,b,c\}\right)$ such that $\alpha(a) = \alpha(b) = \alpha(c) = b$ and the subset of  $n-c$ endomorphisms $\beta$ from the left part of $\mathcal{L}^{a}_{a}\left(\triangle^{(n)}\{a,b,c\}\right)$ such that $\beta(a) = \beta(b) = \beta(c) = b$. The order of this semiring is  $\left|N^{[b]}\left( \triangle^{(n)}\{a,b,c\}\right)\right| = (a+1)(n-c)$.

  Finally, the semiring $N^{[c]}\left( \triangle^{(n)}\{a,b,c\}\right)$ can be represented as a geometric trapezoid whose ``vertices'' are endomorphisms $b_bc_{n-b}$, $\overline{c}$, $a_{a}c_{n-a}$ and $a_ab_{b-a}c_{n-b}$ and whose ``sides'' are semiring  $N^{[c]}\left(\mathcal{STR}^{(n)}\{b,c\}\right)$,  semiring $N^{[c]}\left(\mathcal{STR}^{(n)}\{a,c\}\right)$,  the subsetset of  $b-a+1$ endomorphisms $\alpha$ from the right part of $\mathcal{L}^{a}_{a}\left(\triangle^{(n)}\{a,b,c\}\right)$ such that $\alpha(b) = \alpha(c) = c$ and the subset of  $a+1$ endomorphisms $\beta$ from the right part of $\mathcal{L}^{n-b}_{c}\left(\triangle^{(n)}\{a,b,c\}\right)$ such that $\beta(a) = b, \beta(b) = \beta(c) = c$. The order of this semiring is  $\left|N^{[c]}\left( \triangle^{(n)}\{a,b,c\}\right)\right| = \frac{1}{2}(a+1)(2b-a+2)$.

\vspace{3mm}

\textbf{Proposition} \z  \textsl{The semirings $N^{[a]}\left( \triangle^{(n)}\{a,b,c\}\right)$, $N^{[b]}\left( \triangle^{(n)}\{a,b,c\}\right)$ and $N^{[c]}\left( \triangle^{(n)}\{a,b,c\}\right)$ are trivial.}

\emph{Proof.} Let $\alpha \in N^{[a]}\left( \triangle^{(n)}\{a,b,c\}\right)$. Then $a$ is a unique fixed point of $\alpha$. If we assume that $\alpha(b) = c$, then $c \geq \alpha(c) \geq \alpha(b) =c$ implies $\alpha(c) = c$, which is impossible. So, $\alpha(b) = a$. Now $\alpha(c) = b$, or $\alpha(c) = a$, which means that $\alpha = \overline{a}$. For any $\beta \in N^{[a]}\left( \triangle^{(n)}\{a,b,c\}\right)$ we find $(\alpha\cdot\beta)(a) = a$, $(\alpha\cdot\beta)(b) = \beta(\alpha(b)) = \beta(a) = a$ and $(\alpha\cdot\beta)(c) = \beta(\alpha(c)) = \beta(b) = a$. Hence, $\alpha\cdot \beta = \overline{a}$ and $N^{[a]}\left( \triangle^{(n)}\{a,b,c\}\right)$ is a trivial semiring. The same argument shows that $N^{[b]}\left( \triangle^{(n)}\{a,b,c\}\right)$ and $N^{[c]}\left( \triangle^{(n)}\{a,b,c\}\right)$ are also trivial semirings. \hfill $\Box$

\vspace{3mm}

\emph{{Example}} \z  Let us consider  triangle $\triangle^{(6)}\{1,3,4\}$. Figure 5 illustrates the semirings of 1--nilpotent, 3--nilpotent and 4--nilpotent endomorphisms, an idempotent triangle and right identities. Here we observe that $\triangle^{(6)}\{1,3,4\}$ can be represented as a union of the following semirings: $N^{[1]}\left( \triangle^{(6)}\{1,3,4\}\right)$, $N^{[3]}\left( \triangle^{(6)}\{1,3,4\}\right)$, $N^{[4]}\left( \triangle^{(6)}\{1,3,4\}\right)$, $\mathcal{L}^{2}_{1}\left(\triangle^{(6)}\{1,3,4\}\right)$, $\mathcal{L}^{3}_{1}\left(\triangle^{(6)}\{1,3,4\}\right)$ and $\mathcal{L}^{2}_{4}\left(\triangle^{(6)}\{1,3,4\}\right)$. But this union is not disjoint since the right identities are intersections of the basic layers of the triangle.

\begin{figure}[h]\centering
  \includegraphics[width=105mm]{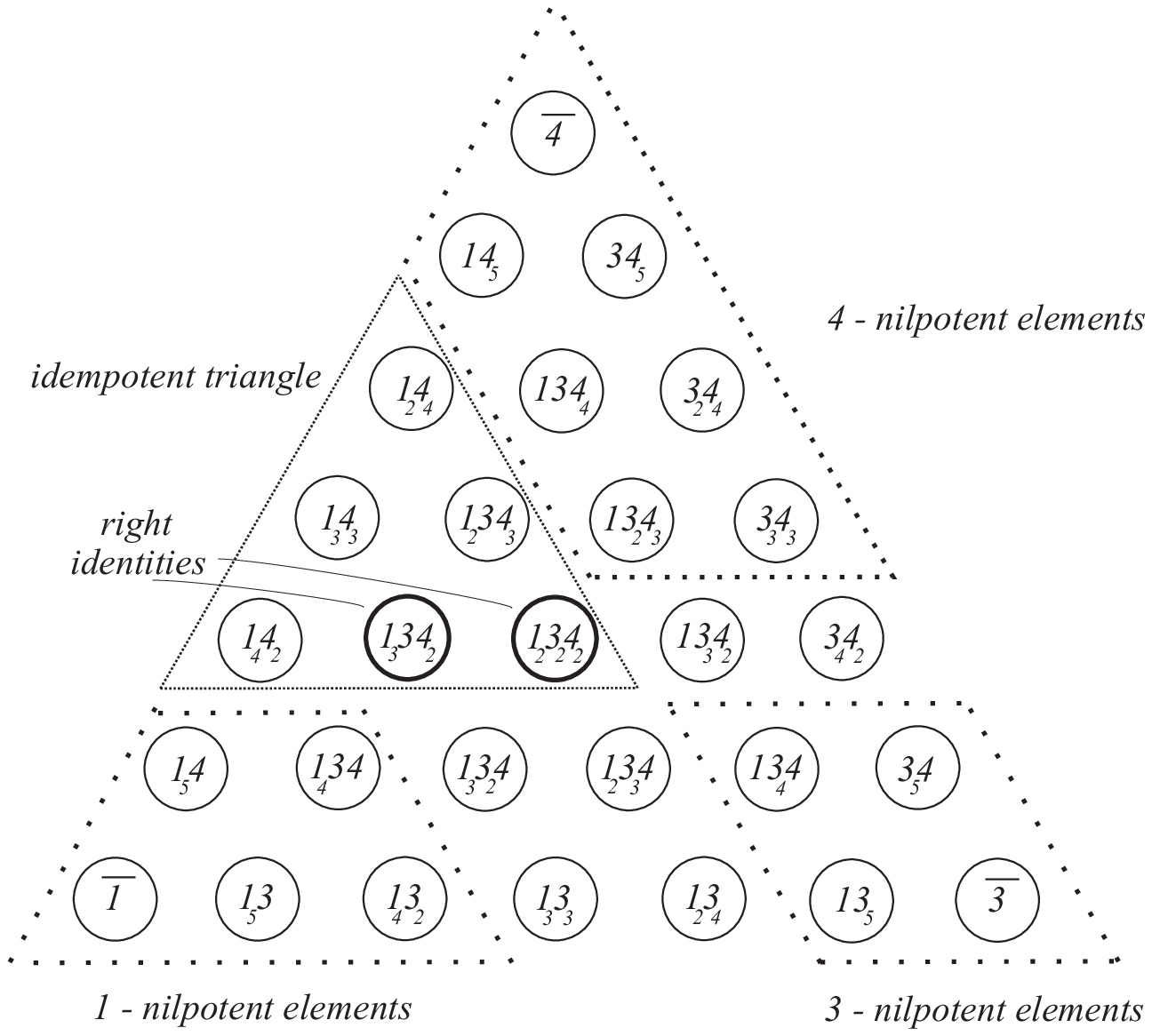}\\
\end{figure}

\vspace{-1mm}

\centerline{\small Figure 5.}
\vspace{4mm}

In order to represent $\triangle^{(n)}\{a,b,c\}$ as a disjoint union of its subsemirings, we look for new subsemirings. Now we consider the set of all the left elements of the basic layers with respect to $\overline{a}$. This set can be represented as a geometric parallelogram whose ``vertices'' are endomorphisms $a_bb_{n-b}$, $a_{a+1}b_{n-a-1}$, $a_{a+1}b_{c-a}c_{n-c-1}$ and $a_bb_{c-b+1}c_{n-c-1}$. The ``sides'' of this parallelogram are  semiring $Id\left(\mathcal{STR}^{(n)}\{a,b\}\right)$, the set of the left elements of the biggest basic layer $\mathcal{L}^{a+1}_{a}\left(\triangle^{(n)}\{a,b,c\}\right)$, the subset of $b-a$ endomorphisms $\alpha$ of the layer $\mathcal{L}^{n-c-1}_{c}\left(\triangle^{(n)}\{a,b,c\}\right)$ with fixed points $a$ and $b$ and the set of the left elements of the  least basic layer $\mathcal{L}^{b}_{a}\left(\triangle^{(n)}\{a,b,c\}\right)$. We denote this  parallelogram by $L_{{par}}$. Then $\left|L_{{par}}\right| = (b-a)(n-c)$.

 Similarly, we consider the set of all the right elements of the basic layers with respect to $\overline{c}$. This set  can also be represented as a geometric parallelogram whose ``vertices'' are endomorphisms $a_ab_{c-a}c_{n-c}$, $b_{c}c_{n-c}$, $b_{b+1}c_{n-b+1}$ and $a_ab_{b-a+1}c_{n-b-1}$. The ``sides'' of parallelogram are   the set of right elements of the biggest basic layer $\mathcal{L}^{n-c}_{a}\left(\triangle^{(n)}\{a,b,c\}\right)$,
  the semiring $Id\left(\mathcal{STR}^{(n)}\{b,c\}\right)$, the set of the right elements of   least basic layer $\mathcal{L}^{n-b-1}_{c}\left(\triangle^{(n)}\{a,b,c\}\right)$ and
    the subset of $c-b$ endomorphisms $\alpha$ of  layer $\mathcal{L}^{a}_{a}\left(\triangle^{(n)}\{a,b,c\}\right)$ with fixed points $b$ and $c$. We denote this  parallelogram by $R_{{par}}$. Then $\left|R_{{par}}\right| = (a+1)(c-b)$.

\vspace{3mm}

\textbf{Proposition} \z  \textsl{The geometric parallelograms $L_{{par}}$ and $R_{{par}}$ are subsemirings of  triangle $\triangle^{(n)}\{a,b,c\}$.}

\emph{Proof.} We shall prove only that $L_{{par}}$ is a semiring, since the proof for $R_{{par}}$ is the same. Let $\alpha \in L_{\mbox{par}}$. Then $\alpha(a) = a$ and $\alpha(b) = \alpha(c) = b$. It is evident that $(\alpha + \beta)(a) = a$ and $(\alpha + \beta)(b) = (\alpha + \beta)(c) = b$. For products we show  $(\alpha\cdot \beta)(a) = a$, $(\alpha\cdot \beta)(b) = b$ and $(\alpha\cdot \beta)(c) = \beta(\alpha(c)) = \beta(b) = b$. \hfill $\Box$

\vspace{3mm}

As we have seen in the last proofs, any endomorphism $\alpha$ of some subsemiring of the triangle can be characterized by ordered triple $(x,y,z)$, where $\alpha(a) = x$, $\alpha(b) = y$ and $\alpha(c) = z$ and $x, y, z \in \{a,b,c\}$. This triple is called a {\emph{type}} of semiring.

Now we can summarize the results of Theorem 33, corollaries 34 and 35 and propositions 36 and 39 and arrange the following ``puzzle'' -- fig. 6, where we register the type of semirings.

\begin{figure}[h]\centering
  \includegraphics[width=75mm]{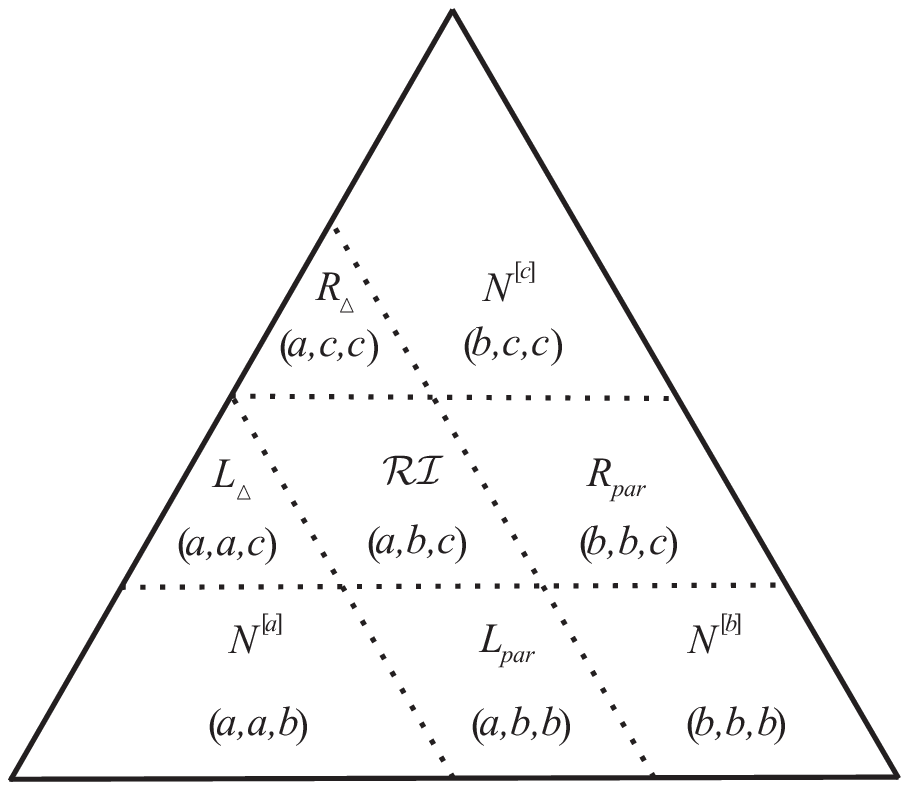}\\
\end{figure}

\vspace{-1mm}

\centerline{\small Figure 6.}
\vspace{4mm}

Actually we prove the following

\vspace{3mm}

\textbf{Theorem} \z \textsl{Any triangle $\triangle^{(n)}\{a,b,c\}$, $n \geq 3$, is a disjoint union of the following subsemirings of the triangle:
 $N^{[a]}\left( \triangle^{(n)}\{a,b,c\}\right)$, $N^{[b]}\left( \triangle^{(n)}\{a,b,c\}\right)$, $N^{[c]}\left( \triangle^{(n)}\{a,b,c\}\right)$, $L_{{par}}$, $R_{{par}}$, $L_{\triangle}$, $R_{\triangle}$ and $\mathcal{RI}\left(\triangle^{(n)}\{a,b,c\}\right)$.}

\vspace{3mm}

As a direct consequence of the last theorem, it follows that  semiring $\triangle^{(n)}\{a,b,c\}\cap {\mathcal{E}}^{(b)}_{\mathcal{C}_n}$ is a disjoint union of semirings $N^{[a]}\left( \triangle^{(n)}\{a,b,c\}\right)$, $L_{{par}}$, $R_{{par}}$ and $\mathcal{RI}\left(\triangle^{(n)}\{a,b,c\}\right)$, that is this semiring can be represented as a geometric parallelogram consisting of the four parallelograms corresponding to these subsemiring -- fig. 6.

\vspace{5mm}

\noindent{\large \bf References}

\vspace{3mm}

[1] J. Je$\hat{\mbox{z}}$ek, T. Kepka and  M. Mar\`{o}ti, ``The endomorphism semiring of a se\-milattice'', \emph{Semigroup Forum}, 78, pp. 21 -- 26, 2009.

[2] I. Trendafilov and D. Vladeva, ``Idempotent Elements of the Endomorphism Semiring of a Finite Chain'', \emph{ISRN Algebra}, vol. 2013, Article ID 120231, 9 pages, 2013.

[3]  I. Trendafilov and  D. Vladeva, ``Nilpotent elements of the endomorphism semiring of a finite chain and Catalan numbers'',
\emph{Proceedings of the Forty Second Spring Conference of the Union of Bulgarian Mathematicians}, Borovetz, April 2--6,  pp. 265 -- 271, 2013.

[4] I. Trendafilov and D. Vladeva, ``Endomorphism semirings without zero of a finite chain'', \emph{Proceedings of the Technical University of Sofia}, vol. 61, no. 2, pp. 9 -- 18, 2011.

[5] I. Trendafilov and D. Vladeva, ``The endomorphism semiring of a finite chain'', \emph{Proceedings of the Technical University of Sofia}, vol. 61, no. 1, pp. 9--18, 2011.

[6] I. Trendafilov and D. Vladeva, ``Subsemirings of the endomorphism semiring of a finite chain'', \emph{Proceedings of the Technical University of Sofia}, vol. 61, no. 1, pp. 19--28, 2011.

[7] J. Zumbr\"{a}gel, ``Classification of finite congruence-simple semirings
with zero,'' \emph{Journal of Algebra and Its Applications}, vol. 7, no. 3, pp. 363--377, 2008.

[8] J.  Golan, \emph{Semirings and Their Applications}, Kluwer, Dordrecht, 1999.

[9] G.  Gratzer, \emph{Lattice Theory: Foundation}, Birkh\"{a}user
Springer Basel AG, 2011.

[10] D. Ferrario and R. Piccinini, \emph{Simplicial Structures in Topology}, Springer, New York, USA, 2011.

[11] M. Desbrun, A, Hirani, M. Leok, J. Marsden, ``Discrete Exterior Calculus'',  arXiv:math/0508341v2 [math.DG] 18 Aug 2005.

\end{document}